# Three-Phase Unbalance Mitigation via DSO-FRA Coordination: A GNB-Based Chance-Constrained Model Considering PV Uncertainty

Qun Zhou, Jiale Guo, Jing Huang*, Minrui Leng, Xueshan Liu, Chenjia Gu*

(College of Electrical Engineering, Sichuan University, Chengdu 610065, China)

* Corresponding authors: huangjing123@scu.edu.cn (Jing Huang), gcj0629@scu.edu.cn (Chenjia Gu)

**Abstract:** Although flexible resource aggregators (FRAs), represented by electric vehicle aggregators (EVAs) and load aggregators (LAs), possess inherent flexibility to mitigate three-phase unbalance, the key challenge lies in how to effectively incentivize them to actively participate in unbalance mitigation. To address this, this paper proposes a generalized Nash bargaining (GNB)-based chance-constrained coordinated operation model for the distribution system operator (DSO) and FRAs. In this framework, the DSO and FRAs cooperate to mitigate three-phase unbalance, where FRAs provide flexibility to reduce the DSO's unbalance mitigation costs, and the DSO offers economic incentives to FRAs in return. Bargaining power is tailored according to each participant's contribution to unbalance mitigation to ensure fair profit allocation. Moreover, a scenario-based chance-constrained formulation is integrated to handle photovoltaic (PV) output uncertainties. To solve the model, a distributed proximal decomposition algorithm (PDA) is introduced for the independent model, while the coordinated model is decomposed into a social welfare maximization subproblem and a payment bargaining subproblem, with an improved bilinear Benders decomposition algorithm developed to solve the former. Numerical results validate the effectiveness of the proposed method in mitigating unbalance, ensuring fair profit distribution, and enhancing robustness against uncertainties.

## 1. Introduction

Three-phase unbalance is one of the most critical issues in power distribution networks (PDNs) [1]. In the UK, more than 70% of low-voltage distribution networks suffer from severe three-phase unbalance [2]. Although PDNs are initially planned with full consideration of potential phase unbalance conditions, the development of PDNs has led to the increasing integration of distributed generation and new-type loads, which has further exacerbated the three-phase unbalance problem [3]. On the one hand, these new elements are not necessarily evenly distributed across the three phases, thereby intensifying the unbalance in the total three-phase load distribution. On the other hand, their high but heterogeneous volatility results in rapid variations of phase unbalance over time, making the unbalance mitigation more challenging. Severe three-phase unbalance can cause numerous adverse effects on PDNs, including increased network losses [4], reduced motor lifespan [5], wasted transformer capacity [6], and malfunction of relay protection devices [7]. It is reported that polyphase induction motors experience additional energy losses exceeding 10% at 2% unbalance, which rise sharply to over 40% at 6% unbalance [8]. Therefore, developing an effective three-phase unbalance mitigation method is of great significance for distribution system operators (DSOs).

Traditional three-phase unbalance mitigation methods mainly fall into three categories: load phase switching

[9], power electronic compensation devices [10] and network reconfiguration [11]. In load phase switching, manual phase switching requires prolonged power outages and struggles to cope with frequent load fluctuations. Automatic phase switching through phase switch devices (PSDs) represents an emerging approach [12]. However, it necessitates substantial investment in PSDs and may compromise power supply reliability. In terms of compensation devices, static var generators (SVGs) and static var compensators (SVCs) are typically installed at transformers to compensate for unbalanced currents [13]. Nevertheless, these power electronic compensators are limited to regulating mild unbalance conditions and cannot cope with severe unbalance scenarios. Moreover, they can only alleviate the unbalance at the transformer level, without addressing the unbalance in the downstream network. Network reconfiguration changes the topology of the PDN by adjusting the open/closed states of tie switches, thereby adjusting the load distribution and mitigating three-phase unbalance [14]. However, network reconfiguration is constrained by the available switching actions and may not be sufficient to address severe or rapidly varying unbalance conditions, and frequent reconfiguration operations may also increase wear on switching devices and compromise system reliability.

Flexible resource aggregators (FRAs), represented by electric vehicle aggregators (EVAs) and load aggregators (LAs), possess inherent flexibility in their consumption patterns. This flexibility can be strategically leveraged to provide phase balancing services. However, current research on FRAs has primarily focused on optimizing their operations to alleviate line congestion [15], shave peak loads [16], eliminate voltage violations [17], and provide ancillary services [18], etc. These studies often assume that the PDN is three-phase balanced, thereby neglecting the impacts of three-phase unbalance. It is noteworthy that even if the three-phase voltages and powers of the PDN remain within the prescribed limits, the degree of unbalance may still far exceed the allowable range, thereby jeopardizing the safe and reliable operation of the PDN.

Some studies have begun to explore three-phase unbalance mitigation by optimizing the operation of FRAs. Reference [19] proposes a data-driven orderly charging framework for electric vehicles (EVs) to mitigate three-phase unbalance and enhance the hosting capacity of photovoltaics (PVs). Reference [20] presents an EV charging strategy that simultaneously addresses both three-phase unbalance and harmonic pollution. Reference [21] establishes a coordination model of EVs and energy storage systems (ESSs) specifically for three-phase unbalance mitigation. Reference [22] further considers the coordinated operation of PVs, ESSs, and LAs for three-phase unbalance mitigation. The above studies verify that strategically adjusting the operation strategies of FRAs can effectively mitigate three-phase unbalance. However, they generally assume that FRAs are able to adjust their operation strategies unconditionally in response to the unbalance mitigation requirements of the PDN. In practice, since the primary objective of FRAs is to minimize their own costs, participating in unbalance mitigation typically incurs additional costs, making the above strategies difficult to apply directly [23]. Therefore, it is necessary to design a reasonable incentive mechanism that aligns the objectives of FRAs with the DSO's unbalance mitigation requirements, thereby motivating FRAs to actively participate in unbalance mitigation.

Few studies have explored incentive mechanisms for FRAs participating in three-phase unbalance mitigation. Reference [24] develops a non-cooperative game-based incentive mechanism that determines compensation for EVs based on the degree of unbalance reduction. Nevertheless, this method primarily focuses on EV gaming on a single bus and fails to consider that EVs may be connected to multiple buses in the PDN, thereby limiting its practical applicability. To address this issue, reference [25] further develops a non-cooperative game model that explicitly considers network topology, incorporating all EV users in the PDN into the gaming framework. However, these existing studies approach FRA participation from a non-cooperative perspective, where individual profit maximization often leads to reduced overall social welfare [26]. In contrast, cooperative game theory has been widely adopted in power system optimization to coordinate the conflicting objectives of multiple self-interested entities. For instance, reference [27] develops a cooperative transactive multi-carrier energy control mechanism for P2P energy and reserve trading under renewable uncertainty based on Nash bargaining theory. Reference [28] proposes a P2P trading mechanism for multi-energy prosumers based on generalized Nash bargaining (GNB), demonstrating that cooperative frameworks can simultaneously enhance individual profits and social welfare. Reference [29] applies Nash bargaining to the alliance expansion of multiple island microgrids with shared hybrid energy storage systems, achieving a win-win outcome through cooperative benefit distribution. Reference [30] further develops an energy-carbon pricing-guided collaborative optimization model for local integrated energy communities, utilizing Nash bargaining theory to coordinate electricity, gas, and heat operations. These studies collectively demonstrate that cooperative game frameworks are capable of reconciling individual rationality with collective optimality, thereby yielding greater overall benefits than non-cooperative approaches [31]. Despite the proven effectiveness of cooperative game theory in general power system coordination, its application to three-phase unbalance mitigation remains largely unexplored. Moreover, most of the above studies distribute the benefits equally among participants, neglecting the varying contributions of different participants to the coalition. This may lead to unfair profit allocation and thereby hinder cooperation. Therefore, from a cooperative perspective, it is of great significance to establish a coordinated operation model between the DSO and FRAs with a fair benefit distribution mechanism.

In addition, the aforementioned studies typically address three-phase unbalance mitigation under a deterministic framework. However, the intermittency and volatility of renewable energy output can cause significant variations in PDN operating states, thereby introducing considerable uncertainty into three-phase unbalance mitigation [32]. Consequently, operation strategies that neglect such variability may lead to economically inefficient or even infeasible solutions. To address this uncertainty, existing studies have predominantly employed robust optimization (RO) [33], distributionally robust optimization (DRO) [34], and stochastic optimization (SO) [35]. For instance, reference [36] establishes a stochastic multi-objective optimization model for PV-integrated three-phase distribution networks, achieving simultaneous minimization of power losses and voltage unbalance under PV uncertainty. Reference [37] proposes a stochastic optimization model for energy management of unbalanced three-phase AC microgrids, ensuring resilient operation under multiple uncertainties. Reference [38] develops a robust three-phase state estimation framework for PV-integrated unbalanced distribution systems, enhancing resilience against bad data. Reference [39] proposes a distributionally robust optimization model for power grids along mountain railways, addressing three-phase unbalance of traction loads under renewable uncertainty. Yet these approaches typically ensure system feasibility through strict checks of all possible or worst-case conditions, then proceed to economic optimization. Considering that the three-phase unbalance degree is not a rigid constraint that must be satisfied at all times, minor and short-term violations of the unbalance limit do not necessarily cause severe impacts on the PDN. However, mitigating such secondary violations can lead to prohibitively high investment costs, making the plans difficult to implement. Therefore, it is necessary to relax the unbalance constraints in some extreme scenarios that occur with low probability, so as to improve the economic efficiency of mitigation without significantly compromising its effectiveness.

Table 1 compares prior studies on three-phase unbalance mitigation in detail and identifies their major differences in mitigation methods and uncertainty characterization. Based on the literature reviewed above, this paper targets the following research gaps:

1) Most existing studies either directly regulate the operation of FRAs by assuming that they will unconditionally follow the DSO's requirements, or adopt non-cooperative game frameworks to incentivize FRA participation. However, direct regulation is impractical because FRAs are self-interested entities that prioritize their own economic objectives, while non-cooperative approaches often lead to conflicts between individual rationality and collective optimality, resulting in reduced overall social welfare. Therefore, it is necessary to establish a cooperative incentive mechanism that aligns the objectives of the DSO and FRAs, enabling them to jointly mitigate three-phase unbalance while achieving mutual benefits.

2) Even when cooperation is considered, the resulting benefits are typically distributed equally among participants or allocated in a simple proportional manner, without accounting for the varying contributions of different participants to unbalance mitigation. This may lead to unfair profit allocation and discourage participants with greater flexibility contributions from active involvement, thereby undermining the stability of the cooperative coalition. Therefore, a fair benefit distribution mechanism that differentiates bargaining power according to each participant's contribution is needed to sustain a stable and effective cooperative relationship.

3) Most existing studies either ignore the impact of renewable energy output uncertainty on unbalance mitigation, thereby leading to infeasible solutions, or handle it using SO, RO, or DRO. However, enforcing feasibility across all possible or worst-case scenarios often gives rise to overly conservative investment strategies, which in turn lead to excessive mitigation costs owing to the high capital expenditure of various assets. Therefore, alternative uncertainty treatment methods are needed, such as chance-constrained programming, which permits three-phase unbalance limit violations at a certain risk level, thereby striking a balance between the economics and effectiveness of unbalance mitigation.

**Table 1** Summary of the reviewed literature for three-phase unbalance mitigation

| Refs | Mitigation method | Incentive mechanism | FRA participation | Social welfare | Profit allocation | Uncertainty handling |
|---|---|---|---|---|---|---|
| [12] | PSD | - | - | - | - | SO |
| [13] | PSD+SVC | - | - | - | - | - |
| [14] | Network reconfiguration | - | - | - | - | - |
| [19]-[22] | FRA optimized operation | Direct regulation | √ | - | - | - |
| [24], [25] | FRA non-cooperative game | Compensation for phase balancing | √ | Non-maximization | Equal allocation | - |
| [36][37] | SVC+ESS | - | - | - | - | SO |
| [39] | ESS | - | - | - | - | DRO |
| Proposed | DSO-FRA coordination | GNB-based | √ | Maximization | Fair allocation | Chance-constrained |

To address these challenges, a GNB-based DSO-FRA chance-constrained coordinated operation model is proposed in this paper for three-phase unbalance mitigation. In this model, the DSO and FRAs cooperate to mitigate three-phase unbalance. FRAs provide flexibility to the DSO to alleviate unbalance at a lower cost. In return, the DSO offers economic incentives to FRAs to strengthen cooperation. To promote fair benefit distribution, a mutually beneficial negotiation agreement is established based on GNB theory. Moreover, PV uncertainties are explicitly considered in the coordinated operation, and a chance-constrained programming formulation is employed to achieve a balance between the robustness and economic efficiency of the coordinated scheme.

The main contributions are summarized as follows:

1) A GNB-based coordinated operation model for the DSO and FRAs is proposed. Through the coordinated operation of the DSO and FRAs, the flexibility of FRAs is fully exploited to effectively mitigate three-phase unbalance while reducing the DSO's investment cost and the FRAs' operation costs, thereby maximizing social welfare. To accurately assess the post-cooperation benefits, independent operation models of the DSO and FRAs are established to determine the disagreement point. Furthermore, the bargaining power of each participant is tailored according to its individual contribution to unbalance mitigation, thereby achieving an equitable profit allocation.

2) To address the impact of PV output uncertainties on cooperation, a scenario-based chance-constrained model is integrated into the proposed DSO-FRA coordinated framework. Specifically, the chance-constrained model does not require the three-phase unbalance constraints to be satisfied in all scenarios; rather, it allows the three-phase unbalance degree to exceed its limit in scenarios with a small probability. By enabling the coordination scheme to effectively mitigate three-phase unbalance under most scenarios, the chance-constrained model significantly enhances the robustness of the coordination scheme and assists the DSO in striking a trade-off between economic efficiency and operational security, thereby facilitating more informed decision-making.

3) Solution algorithms are designed for the independent operation models and the coordinated operation model, respectively. First, A distributed proximal decomposition algorithm (PDA) is introduced to solve the Nash equilibrium problem (NEP) in the FRA independent operation model. Then, based on the decomposable structure of the coordinated operation model, it is decomposed into a social welfare maximization subproblem and a payment bargaining subproblem. Finally, to address the binary variables and bilinear constraints introduced by the chance-constrained model, an improved bilinear Benders decomposition algorithm is proposed to solve the social welfare maximization subproblem.

The rest of this paper is organized as follows. Section 2 presents the problem description. Section 3

formulates the independent operation models of FRAs and the DSO. Section 4 formulates the proposed DSO-FRA coordinated operation model. Section 5 presents the solution methods. Section 6 provides numerical results. Finally, Section 7 concludes the paper.

## 2. Problem description

With the large-scale integration of single-phase distributed generation and new-type loads [40], three-phase unbalance has become one of the most critical power quality issues in PDNs. To address this issue, DSOs typically deploy compensation devices at the planning level [41], among which soft open points (SOPs) have attracted increasing attention in recent years due to their flexible power flow regulation capability. By appropriately configuring the capacity and placement of SOPs, the DSO can effectively transfer active power among phases and provide reactive power support, thereby mitigating three-phase unbalance. However, the investment cost of SOPs remains substantial, and their economic viability is often unsatisfactory. This motivates the exploration of alternative or supplementary measures that can achieve effective unbalance mitigation at a lower cost.

Meanwhile, FRAs, including EVAs and LAs, have emerged as important participants in PDN operations. EVAs primarily manage a large number of single-phase slow-charging EVs, whose charging loads can be shifted across time to provide flexibility. LAs aggregate various types of residential and commercial loads, such as air conditioners, washers, and heaters, which can be adjusted through shiftable and curtailable load management to offer demand-side flexibility. Although these resources are primarily managed for individual economic objectives in independent operation, they possess significant potential to support system-level objectives such as unbalance mitigation when properly coordinated. Fig. 1 illustrates a typical PDN that incorporates FRAs, PVs, and SOPs.

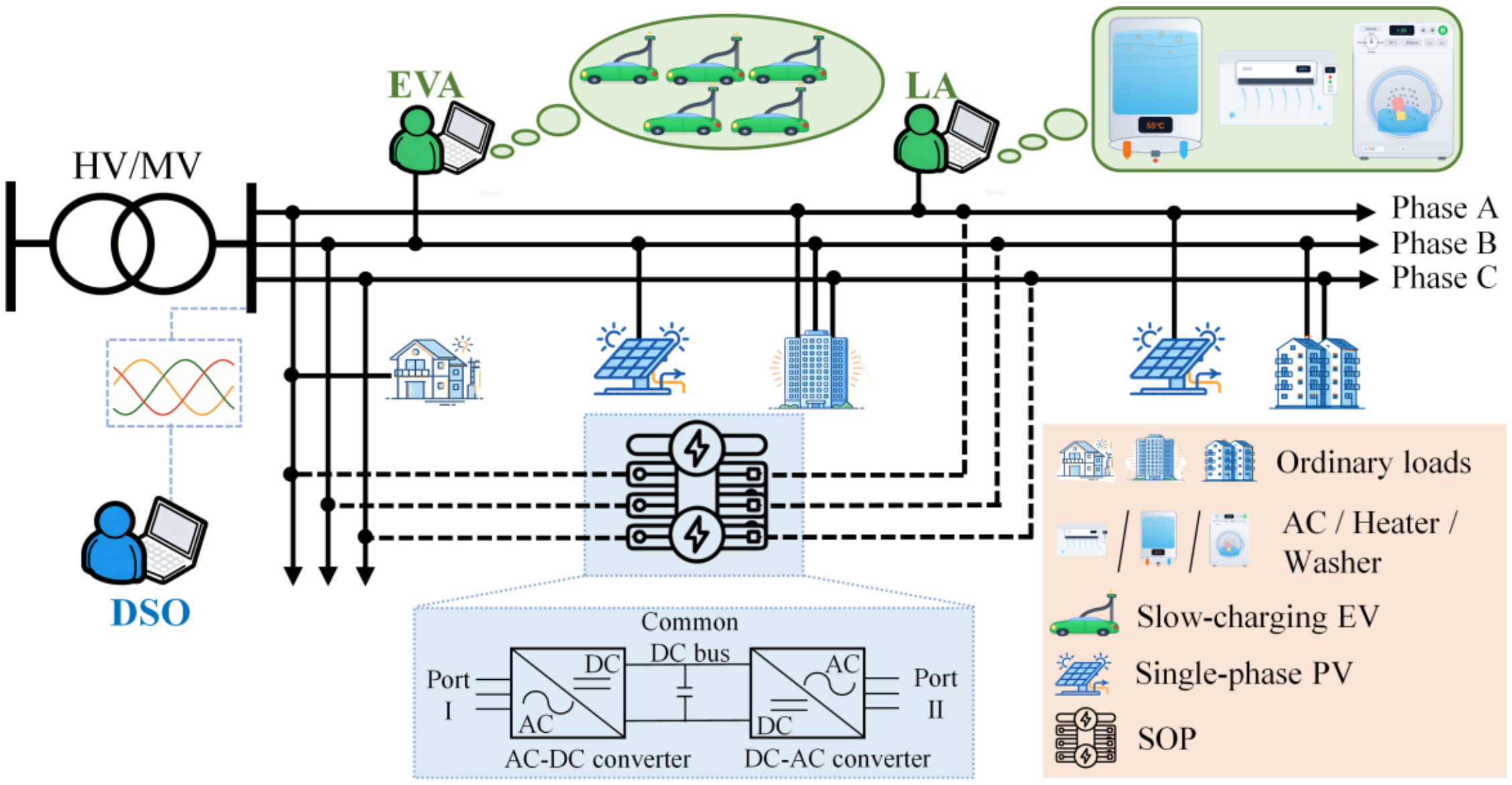

**Fig. 1.** Schematic diagram of a PDN with FRAs, PVs, and SOPs.

To this end, this paper proposes a DSO-FRA coordination framework for three-phase unbalance mitigation, as illustrated in Fig. 2. In the proposed framework, FRAs cooperate with the DSO to mitigate three-phase unbalance by adjusting their consumption patterns, which in turn helps the DSO reduce the investment cost of SOPs. In return, the DSO provides economic incentives to FRAs based on the cost savings achieved through cooperation compared to independent operation. To establish a stable cooperative relationship and ensure a fair distribution of profits, the GNB theory is adopted, allowing for differentiated bargaining power to reflect the varying contributions of different participants to unbalance mitigation.

In summary, the proposed coordination framework involves the following key points. First, it is necessary to establish independent operation models of the DSO and FRAs considering PV output uncertainties, so as to accurately determine the economic benefits of coordination for the DSO and FRAs. Second, the contributions of each participant to unbalance mitigation should be reasonably quantified to achieve a fair profit allocation. Finally, effective solution algorithms need to be developed to efficiently solve this large-scale and complex problem. These aspects are discussed in detail in the following sections.

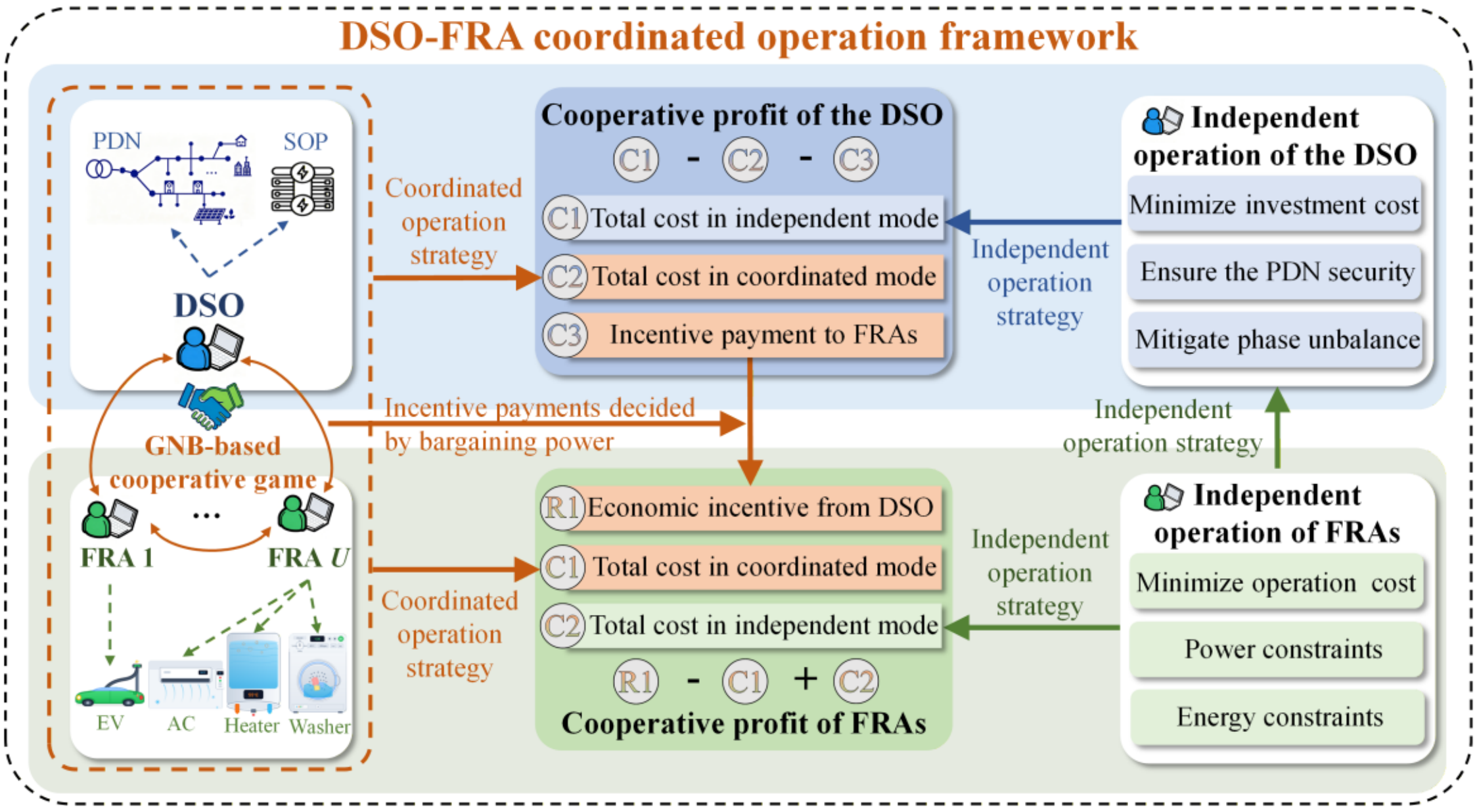

**Fig. 2.** The schematic diagram of the DSO-FRA coordinated operation framework.

## 3. Independent operation models of FRAs and the DSO

### *3.1. Independent operation model of FRAs*

#### *3.1.1. Independent operation model of EVAs*

The optimization objective of each EVA is to minimize its own operation cost, which includes the charging cost of EVs and the cost associated with charging energy deviation. Therefore, the objective function of EVA $n$ can be expressed as follows:

$$\min C_{n,s}^{\text{EVA}} = C_{n,s}^{\text{cha}} + C_{n,s}^{\text{dev}} \tag{1}$$

$$C_{n,s}^{\text{cha}} = \sum_{t\in T} \left( a_t \left( \sum_{n\in\mathcal{N}} x_{n,t,s} + \sum_{m\in\mathcal{M}} x_{m,t,s} \right) + b_t \right) x_{n,t,s} \Delta t \quad \forall s,n \tag{2}$$

$$C_{n,s}^{\text{dev}} = c_1 e_{n,s}^{\text{dev}} \quad \forall s,n \tag{3}$$

where $C_{n,s}^{\text{cha}}$ is the charging cost of EVA $n$ in scenario $s$. $C_{n,s}^{\text{dev}}$ is the energy deviation cost of EVA $n$ in scenario $s$. $a_t$ and $b_t$ are the dynamic pricing coefficients. $x_{n,t,s}$ is the charging power of EVA $n$ at time $t$ in scenario $s$. $x_{m,t,s}$ is the power of LA $m$ at time $t$ in scenario $s$. $\Delta t$ is the time interval. $c_1$ is the cost coefficient of energy deviation. $e_{n,s}^{\text{dev}}$ is the energy deviation of EVA $n$ in scenario $s$.

Since each EVA manages a large number of EVs, its operational constraints, including charging power limits, charging energy limits, and energy deviation limits, are obtained by aggregating the relevant constraints of individual EVs. These operational constraints of EVA $n$ can be expressed as follows:

$$x_{n,t,s} \le x_n^+ \quad \forall s,t,n \tag{4}$$

$$0 \le x_{n,t,s} \le \sum_{i\in\mathcal{I}} p_{i,t}^+ \quad \forall s,t,n \tag{5}$$

$$\sum_{i\in\mathcal{I}} e_{i,t}^- \le \sum_{\tau\le t} \eta x_{n,\tau,s}\Delta t \le \sum_{i\in\mathcal{I}} e_{i,t}^+ \quad \forall s,t,n \tag{6}$$

$$e_{n,s}^{\text{dev}} + \sum_{i\in\mathcal{I}} \eta x_{n,t,s}\Delta t = \sum_{i\in\mathcal{I}} \left(e_i^f - e_i^s\right) \quad \forall s,t,n \tag{7}$$

where $x_n^+$ is the maximum charging power of EVA $n$ due to charging station capacity. $p_i^+$ is the maximum charging power of EV $i$ due to battery technology. $\eta$ is the charging efficiency. $e_{i,t}^-$ and $e_{i,t}^+$ are the minimum and maximum charging energy of EV $i$ at time $t$. $e_i^s$ is the initial energy of EV $i$. $e_i^f$ is the expected energy of EV $i$. Constraints (4)–(7) constitute the aggregated charging model of EVA $n$. Specifically, constraints (4)–(5) limit the aggregated charging power of EVA $n$ at each time slot, constraint (6) restricts the aggregated charging energy, and the aggregated energy deviation of EVA $n$ is calculated by (7).

The boundaries of charging power and energy for a single EV $i$ within EVA $n$ are further defined by following constraints:

$$p_{i,t}^+ = \begin{cases} p_i^+ & \forall t \in \left\{t_i^s,\dots,t_i^f\right\}, \forall i \\ 0 & \forall t \notin \left\{t_i^s,\dots,t_i^f\right\}, \forall i \end{cases} \tag{8}$$

$$e_{i,t}^- = \max\left\{e_i^f - e_i^s - \eta p_{i,t}^+\left(t_i^f - t\right)\Delta t - e_i^{\text{dev}+}, 0\right\} \quad \forall t,i \tag{9}$$

$$e_{i,t}^+ = \min\left\{\eta p_{i,t}^+\left(t - t_i^s + 1\right)\Delta t, e_i^f - e_i^s\right\} \quad \forall t,i \tag{10}$$

where $t_i^s$ and $t_i^f$ are the arrival and departure times of EV $i$. $e_i^{\text{dev}+}$ is the maximum allowable energy deviation for EV $i$.

*3.1.2. Independent operation model of LAs*

Analogous to that of EVAs, the optimization objective of each LA is to minimize its own operation cost, which consists of minimizing the electricity procurement cost of loads and the load curtailment cost. Therefore, the objective function of LA $m$ is formulated as follows:

$$\min C_{m,s}^{\text{LA}} = C_{m,s}^{\text{cha}} + C_{m,s}^{\text{dev}} \tag{11}$$

$$C_{m,s}^{\text{cha}} = \sum_{t\in\mathcal{T}} \left( a_t \left( \sum_{n\in\mathcal{N}} x_{n,t,s} + \sum_{m\in\mathcal{M}} x_{m,t,s} \right) + b_t \right) x_{m,t,s}\Delta t \quad \forall s,m \tag{12}$$

$$C_{m,s}^{\text{dev}} = c_2 e_{m,s}^{\text{dev}} \quad \forall s,m \tag{13}$$

where $C_{m,s}^{\text{cha}}$ is the electricity cost of LA $m$ in scenario $s$. $C_{m,s}^{\text{dev}}$ is the load curtailment cost of LA $m$ in scenario $s$. $c_2$ is the load curtailment cost coefficient. $e_{m,s}^{\mathrm{dev}}$ is the load curtailment energy of LA $m$ in scenario $s$.

The loads managed by each LA include base load, shiftable load, and curtailable load. Among them, the base load is non-dispatchable; the power of the shiftable load can vary, but its total energy consumption over a certain period must remain unchanged; and the curtailment of the curtailable load shall not exceed its upper limit. These constraints are formulated as follows:

$$\sum_{t\in\mathcal{T}} p_{m,t,s}^{sh}\Delta t = \sum_{t\in\mathcal{T}} P_{m,t}^{sh}\Delta t \quad \forall s,m \tag{14}$$

$$p_{m,t,s}^{sh} \geq 0 \quad \forall s,t,m \tag{15}$$

$$0 \leq p_{m,t,s}^{cut} \leq P_{m,t}^{cut} \quad \forall s,t,m \tag{16}$$

$$x_{m,t,s} = P_{m,t}^{base} + p_{m,t,s}^{sh} + p_{m,t,s}^{cut} \quad \forall s,t,m \tag{17}$$

$$x_m^{+} = \max_{t\in\mathcal{T}}\left(P_{m,t}^{sh} + P_{m,t}^{cut}\right) \quad \forall m \tag{18}$$

$$e_{m,s}^{\mathrm{dev}} = \sum_{t\in\mathcal{T}} x_{m,t,s}\Delta t - \sum_{t\in\mathcal{T}}\left(P_{m,t}^{base} + P_{m,t}^{shift} + P_{m,t}^{cut}\right)\Delta t \quad \forall s,m \tag{19}$$

where $p_{m,t,s}^{sh}$ is the power of the shiftable load of LA $m$ in scenario $s$ at time $t$. $P_{m,t}^{sh}$ i s the scheduled power of the shiftable load of LA $m$ at time $t$. $p_{m,t,s}^{cut}$ is the power of the curtailable load of LA $m$ in scenario $s$ at time $t$. $P_{m,t}^{cut}$ is the scheduled power of the curtailable load of LA $m$ at time $t$. $P_{m,t}^{base}$ is the base load of LA $m$ at time $t$. $x_m^{+}$ is the maximum schedulable power of LA $m$. Constraints (14)–(19) constitute the scheduling constraints of LA $m$. Specifically, constraint (14) ensures that the total energy consumption of the shiftable loads remains unchanged, constraints (15)–(16) define the power limits of the shiftable and curtailable loads, constraint (17) calculates the aggregated power of LA $m$, constraint (18) defines the maximum schedulable power of LA $m$, which is utilized for the bargaining power calculation, and constraint (19) calculates the load curtailment energy of LA $m$.

*3.1.3. Compact form of independent operation models of FRAs*

The independent operation models of EVAs and LAs are similar in form: the first term is the purchase cost (charging cost for EVAs and electricity cost for LAs), and the second term is the deviation cost (dissatisfaction cost for EVAs due to unmet charging energy, and load curtailment cost for LAs). Therefore, they are unified as the independent operation model **FRA⁰**, which can be expressed as follows:

$$\begin{cases} \min\limits_{\boldsymbol{x}_{u,s}} C_{u,s}^{\mathrm{FRA}}\left(\boldsymbol{x}_{u,s}, \boldsymbol{x}_{-u,s}\right) & \forall s,u \\ \text{s.t.} \quad \boldsymbol{x}_{u,s} \in \mathcal{X}_{n,s}^{\mathrm{EVA}} = \left\{(2)-(10)\right\} & \forall s, \forall u \in \mathcal{U}\cap\mathcal{N} \\ \quad\quad\ \boldsymbol{x}_{u,s} \in \mathcal{X}_{m,s}^{\mathrm{LA}} = \left\{(12)-(19)\right\} & \forall s, \forall u \in \mathcal{U}\cap\mathcal{M} \end{cases} \tag{20}$$

The optimization objective of each FRA $u$ is related not only to its own strategy $\boldsymbol{x}_{u,s} = \left\{x_{u,1,s};\ldots;x_{u,T,s},\forall s\right\}$, but also to the strategies of other FRAs $\boldsymbol{x}_{-u,s} = \left\{\boldsymbol{x}_{1,s};\ldots;\boldsymbol{x}_{u-1,s};\boldsymbol{x}_{u+1,s};\ldots;x_{U,s},\forall s\right\}$. Therefore, the problem **FRA⁰** constitutes a NEP. The existence and uniqueness of the solution have been proven in [42].

*3.2. Independent operation model of the DSO*

*3.2.1 Multiphase power flow model*

In practice, PDNs typically consist of a large number of single-phase loads, unbalanced three-phase loads, and mixed wye- and delta-connected configurations, and may also operate under conditions with missing phases. Therefore, a multiphase power flow model is adopted in this paper to enhance the applicability of the proposed approach to real-world distribution systems.

Let $\dot{\boldsymbol{s}}_j^Y=[\dot{s}_j^A;\dot{s}_j^B;\dot{s}_j^C]$ and $\dot{\boldsymbol{s}}_j^\Delta=[\dot{s}_j^{AB};\dot{s}_j^{BC};\dot{s}_j^{CA}]$ denote the wye- and delta-connected complex power injections at PQ bus $j$, respectively. Let $\dot{\boldsymbol{v}}_j=[\dot{v}_j^A;\dot{v}_j^B;\dot{v}_j^C]$ denote the complex voltage at PQ bus $j$. Let $\dot{\boldsymbol{s}}^Y=[\dot{\boldsymbol{s}}_j^Y]_{j=1}^J$, $\dot{\boldsymbol{s}}^\Delta=[\dot{\boldsymbol{s}}_j^\Delta]_{j=1}^J$, $\dot{\boldsymbol{v}}=[\dot{\boldsymbol{v}}_j]_{j=1}^J$ be the vectors collecting the relevant quantities of all PQ buses in the PDN. Let $\dot{\boldsymbol{s}}_0=[\dot{s}_0^A;\dot{s}_0^B;\dot{s}_0^C]$ and $\dot{\boldsymbol{v}}_0=[\dot{v}_0^A;\dot{v}_0^B;\dot{v}_0^C]$ denote the complex power injection and nodal voltage at the slack bus. Then the multiphase power flow equation of the PDN can be expressed as:

$$\dot{\boldsymbol{s}}_0=\mathrm{diag}\left(\dot{\boldsymbol{v}}_0\right)\left(\bar{\dot{\mathbf{Y}}}_{00}\bar{\dot{\boldsymbol{v}}}_0+\bar{\dot{\mathbf{Y}}}_{0L}\bar{\dot{\boldsymbol{v}}}\right) \tag{21}$$

$$\dot{\boldsymbol{v}}=\dot{\boldsymbol{v}}_0'+\dot{\mathbf{Y}}_{LL}^{-1}\left(\mathrm{diag}\left(\bar{\dot{\boldsymbol{v}}}\right)^{-1}\bar{\dot{\boldsymbol{s}}}^Y+\mathbf{H}^{\mathrm{T}}\mathrm{diag}\left(\mathbf{H}\bar{\dot{\boldsymbol{v}}}\right)^{-1}\bar{\dot{\boldsymbol{s}}}^\Delta\right) \tag{22}$$

where $\dot{\mathbf{Y}}=[\dot{\mathbf{Y}}_{00},\ \dot{\mathbf{Y}}_{0L};\ \dot{\mathbf{Y}}_{L0},\ \dot{\mathbf{Y}}_{LL}]\in\mathbb{C}^{(3+B)\times(3+B)}$ is the three-phase node admittance matrix. $B$ is the total number of phases for all PQ buses. $\dot{\boldsymbol{v}}_0'=-\dot{\mathbf{Y}}_{LL}^{-1}\dot{\mathbf{Y}}_{L0}\dot{\boldsymbol{v}}_0$ is the zero-load voltage vector. $\mathbf{H}\in\mathbb{C}^{L\times B}$ is a block matrix defined by (23). $L$ is the total number of phase-to-phase connections for all PQ buses. In addition, for a PQ bus $j$ with missing phases, the vectors $\dot{\boldsymbol{s}}_j^Y$, $\dot{\boldsymbol{s}}_j^\Delta$, $\dot{\boldsymbol{v}}_j$, and $\boldsymbol{\Gamma}_j$ only contain the quantities corresponding to the existing phases.

$$\mathbf{H}=\begin{bmatrix}\boldsymbol{\Gamma}_1 & & \\ & \ddots & \\ & & \boldsymbol{\Gamma}_J\end{bmatrix}\quad \boldsymbol{\Gamma}_j=\begin{bmatrix}1 & -1 & 0\\ 0 & 1 & -1\\ -1 & 0 & 1\end{bmatrix} \tag{23}$$

It can be observed that (22) is a fixed-point equation with respect to the voltage vector. Assume that $\dot{\boldsymbol{v}}'$ is a set of solutions obtained by solving (22). Then, by taking $\dot{\boldsymbol{v}}'$ as the initial point and performing one iteration of (22), the linearized form of the voltage phasors can be derived as follows:

$$\dot{\boldsymbol{v}}=\dot{\boldsymbol{v}}_0'+\dot{\boldsymbol{K}}^Y\left[\boldsymbol{p}^Y;\boldsymbol{q}^Y\right]+\dot{\boldsymbol{K}}^\Delta\left[\boldsymbol{p}^\Delta;\boldsymbol{q}^\Delta\right] \tag{24}$$

where $\dot{\boldsymbol{K}}^Y=[\dot{\mathbf{Y}}_{LL}^{-1}\mathrm{diag}\left(\bar{\dot{\boldsymbol{v}}}'\right)^{-1},-j\dot{\mathbf{Y}}_{LL}^{-1}\mathrm{diag}\left(\bar{\dot{\boldsymbol{v}}}'\right)^{-1}]$, $\dot{\boldsymbol{K}}^\Delta=[\dot{\mathbf{Y}}_{LL}^{-1}\mathbf{H}^{\mathrm{T}}\mathrm{diag}\left(\mathbf{H}\bar{\dot{\boldsymbol{v}}}'\right)^{-1},-j\dot{\mathbf{Y}}_{LL}^{-1}\mathbf{H}^{\mathrm{T}}\mathrm{diag}\left(\mathbf{H}\bar{\dot{\boldsymbol{v}}}'\right)^{-1}]$, $[\boldsymbol{p}^Y;\boldsymbol{q}^Y]=[\mathrm{Re}\{\dot{\boldsymbol{s}}^Y\};\mathrm{Im}\{\dot{\boldsymbol{s}}^Y\}]$, and $[\boldsymbol{p}^\Delta;\boldsymbol{q}^\Delta]=[\mathrm{Re}\{\dot{\boldsymbol{s}}^\Delta\};\mathrm{Im}\{\dot{\boldsymbol{s}}^\Delta\}]$.

Taking the modulus of both sides of (24) and assuming that the voltage drops are much smaller than the nominal voltage [43], the linearized expression for the voltage amplitude is obtained as follows:

$$\boldsymbol{v}=\left|\dot{\boldsymbol{v}}_0'\right|+\boldsymbol{E}^Y\left[\boldsymbol{p}^Y;\boldsymbol{q}^Y\right]+\boldsymbol{E}^\Delta\left[\boldsymbol{p}^\Delta;\boldsymbol{q}^\Delta\right] \tag{25}$$

where $\boldsymbol{E}^Y=\left|\mathrm{diag}(\dot{\boldsymbol{v}}_0')\right|\mathrm{Re}\{\mathrm{diag}(\dot{\boldsymbol{v}}_0')^{-1}\dot{\boldsymbol{K}}^Y\}$ and $\boldsymbol{E}^\Delta=\left|\mathrm{diag}(\dot{\boldsymbol{v}}_0')\right|\mathrm{Re}\{\mathrm{diag}(\dot{\boldsymbol{v}}_0')^{-1}\dot{\boldsymbol{K}}^\Delta\}$.

In addition, the linear expression for injected active power at slack bus and PDN power loss can be derived as follows:

$$\boldsymbol{p}_0 = \mathrm{Re}\{\dot{\boldsymbol{s}}_0\} = \boldsymbol{p}_0' + \boldsymbol{F}^Y\left[\boldsymbol{p}^Y;\boldsymbol{q}^Y\right] + \boldsymbol{F}^\Delta\left[\boldsymbol{p}^\Delta;\boldsymbol{q}^\Delta\right] \tag{26}$$

$$p^{\mathrm{loss}} = \mathbf{1}_3^{\mathrm{T}}\boldsymbol{p}_0 + \mathbf{1}_B^{\mathrm{T}}\boldsymbol{p}^Y + \mathbf{1}_L^{\mathrm{T}}\boldsymbol{p}^\Delta \tag{27}$$

where $\boldsymbol{p}_0' = \mathrm{Re}\{\mathrm{diag}(\dot{\boldsymbol{v}}_0)(\bar{\dot{\mathbf{Y}}}_{00} - \bar{\dot{\mathbf{Y}}}_{0L}\bar{\dot{\mathbf{Y}}}_{LL}^{-1}\bar{\dot{\mathbf{Y}}}_{L0})\bar{\dot{\boldsymbol{v}}}_0\}$, $\boldsymbol{F}^Y = \mathrm{Re}\{\mathrm{diag}(\dot{\boldsymbol{v}}_0)\bar{\dot{\mathbf{Y}}}_{0L}\bar{\dot{\boldsymbol{K}}}^Y\}$, $\boldsymbol{F}^\Delta = \mathrm{Re}\{\mathrm{diag}(\dot{\boldsymbol{v}}_0)\bar{\dot{\mathbf{Y}}}_{0L}\bar{\dot{\boldsymbol{K}}}^\Delta\}$.

*3.2.2 Chance-constrained independent operation model of the DSO*

The optimization objective of the DSO is to minimize the investment cost of SOPs and the energy loss cost of the PDN, given the operation strategies of FRAs and subject to the operational requirements of the PDN. Therefore, the objective function of the DSO can be expressed as follows:

$$\min C^{\mathrm{DSO}} = C^{\mathrm{inv}} + \sum_{s\in\mathcal{S}} \pi_s C_s^{\mathrm{loss}} \tag{28}$$

$$C^{\mathrm{inv}} = c_3 \frac{1}{365} \frac{r(1+r)^y}{(1+r)^y - 1} \sum_{d\in\mathcal{D}} S_d \tag{29}$$

$$C_s^{\mathrm{loss}} = c_4 \sum_{t\in\mathcal{T}} p_{t,s}^{\mathrm{loss}} \Delta t \quad \forall s \tag{30}$$

where $C^{\mathrm{inv}}$ is the investment cost of SOPs. $C_s^{\mathrm{loss}}$ is the energy loss cost in scenario *s*. $\pi_s$ is the probability of scenario *s*. $c_3$ is the unit investment cost of SOPs per capacity. *r* is the discount rate. *y* is the planning horizon of SOPs. $S_d$ is the planned capacity of SOP *d*. $c_4$ is the energy loss cost coefficient. $p_{t,s}^{\mathrm{loss}}$ is the power loss at time *t* in scenario *s*.

As power electronic devices for power transmission between nodes, SOPs are subject to power balance constraints and capacity constraints during operation. In addition, the maximum installation capacity of SOPs should also be limited. These constraints can be expressed as follows:

$$\sum_{\varphi\in\phi} p_{I,d,\varphi,t,s} + \sum_{\varphi\in\phi} p_{II,d,\varphi,t,s} = 0 \quad \forall s,t,d \tag{31}$$

$$q^- \le q_{I,d,\varphi,t,s} \le q^+ \quad \forall s,t,d,\varphi \tag{32}$$

$$q^- \le q_{II,d,\varphi,t,s} \le q^+ \quad \forall s,t,d,\varphi \tag{33}$$

$$3\sqrt{p_{I,d,\varphi,t,s}^2 + q_{I,d,\varphi,t,s}^2} \le S_d \quad \forall s,t,d,\varphi \tag{34}$$

$$3\sqrt{p_{II,d,\varphi,t,s}^2 + q_{II,d,\varphi,t,s}^2} \le S_d \quad \forall s,t,d,\varphi \tag{35}$$

$$S_d^{\mathrm{SOP}} \le S_d^+ \tag{36}$$

where $p_{I,d,\varphi,t,s}$ and $q_{I,d,\varphi,t,s}$ are the active and reactive power on the *I*-side of SOP *d* in scenario *s* at time *t* for phase $\varphi$. $p_{II,d,\varphi,t,s}$ and $q_{II,d,\varphi,t,s}$ are the active and reactive power on the *II*-side of SOP *d* in scenario *s* at time *t* for phase $\varphi$. $q^+$ and $q^-$ are the upper and lower limits of reactive power output of SOPs. $S_d^+$ is the upper limit of the installation capacity of SOP *d*. Constraints (31)–(36) constitute the planning and operational constraints of SOPs. Specifically, constraint (31) defines the active power flow on both sides of SOPs, constraints (32)–(33) limit the reactive power outputs on both sides of SOPs, constraints (34)–(35) limit the total capacity on both sides of SOPs, and constraint (36) limits the upper bound of SOP planning capacity.

During the process of three-phase unbalance mitigation, the basic operational constraints of the PDN must be satisfied, including nodal power constraints and upper and lower voltage limits. These constraints can be expressed as follows:

$$p_{j,\varphi,t,s} = p_{j,\varphi,t,s}^{base} + \sum_{n\in\mathcal{N}} D_{j,n,\varphi}^{\mathrm{EVA}} x_{n,t,s} + \sum_{m\in\mathcal{M}} D_{j,m,\varphi}^{\mathrm{LA}} x_{m,t,s} + \sum_{d\in\mathcal{D}} D_{j,d}^{\mathrm{SOP}} p_{d,\varphi,t,s} - \sum_{g\in\mathcal{G}} D_{j,g,\varphi}^{\mathrm{PV}} p_{g,t,s} \quad \forall s,t,\varphi,j \tag{37}$$

$$q_{j,\varphi,t,s} = q_{j,\varphi,t,s}^{base} + \xi^{\mathrm{EVA}} \sum_{n\in\mathcal{N}} D_{j,n,\varphi}^{\mathrm{EVA}} x_{n,t,s} + \xi^{\mathrm{LA}} \sum_{m\in\mathcal{M}} D_{j,m,\varphi}^{\mathrm{LA}} x_{m,t,s} + \sum_{d\in\mathcal{D}} D_{j,d}^{\mathrm{SOP}} q_{d,\varphi,t,s} - \xi^{\mathrm{PV}} \sum_{g\in\mathcal{G}} D_{j,g,\varphi}^{\mathrm{PV}} p_{g,t,s} \quad \forall s,t,\varphi,j \tag{38}$$

$$p_{j,\varphi,t,s}^{base} + \sum_{n\in\mathcal{N}} D_{j,n,\varphi}^{\mathrm{EVA}} x_{n,t,s} + \sum_{m\in\mathcal{M}} D_{j,m,\varphi}^{\mathrm{LA}} x_{m,t,s} + \sum_{d\in\mathcal{D}} D_{j,d}^{\mathrm{SOP}} p_{d,\varphi,t,s} \ge 0 \quad \forall s,t,\varphi,j \tag{39}$$

$$v^{-} \le v_{j,\varphi,t,s} \le v^{+} \quad \forall s,t,\varphi,j \tag{40}$$

where $p_{j,\varphi,t,s}$ and $q_{j,\varphi,t,s}$ are the total active and reactive power at bus $j$ in scenario $s$ at time $t$ for phase $\varphi$. $p_{j,\varphi,t,s}^{base}$ and $q_{j,\varphi,t,s}^{base}$ are the base load active and reactive power at bus $j$ in scenario $s$ at time $t$ for phase $\varphi$. $p_{g,t,s}$ is the active power output of PV $g$ in scenario $s$ at time $t$. $D_{j,n,\varphi}^{\mathrm{EVA}}$, $D_{j,m,\varphi}^{\mathrm{LA}}$, and $D_{j,g,\varphi}^{\mathrm{PV}}$ are binary constants equal to 1 if EVA $n$, LA $m$, and PV $g$ are connected to bus $j$ in phase $\varphi$. $D_{j,d}^{\mathrm{SOP}}$ is a binary constant equal to 1 if SOP $d$ is located at bus $j$. $\xi^{\mathrm{EVA}}$, $\xi^{\mathrm{LA}}$, and $\xi^{\mathrm{PV}}$ are coefficients for calculating reactive power of EVAs, LAs, and PVs. $v^{+}$ and $v^{-}$ are the upper and lower limits of bus voltages. $v_{j,\varphi,t,s}$ is the voltage magnitude of phase $\varphi$ of bus $j$ in scenario $s$ at time $t$. Constraints (37)–(40) constitute the network operational constraints. Specifically, constraints (37)–(38) calculate the nodal loads of the PDN, constraint (39) limits the active power that can be transmitted by SOPs, constraint (40) limits the nodal voltages of the PDN.

To avoid over-investment in SOPs, a chance-constrained formulation is further introduced in this paper. By allowing a limited probability of unbalance violations in low-probability extreme scenarios, while ensuring compliance under most operating conditions, the chance constraint provides a flexible trade-off between mitigation effectiveness and economic efficiency. The chance constraint is formulated as follows:

$$\Pr\left\{U_{j,t,s} \le U^{+}\right\} \ge 1-\varepsilon \quad \forall s,t,j \tag{41}$$

$$U_{j,t,s} = \max_{\varphi\in\phi} \frac{\left|v_{j,\varphi,t,s} - v_{j,t,s}^{avg}\right|}{v_{j,t,s}^{avg}} \quad \forall s,t,j \tag{42}$$

where $U_{j,t,s}$ is the three-phase unbalance degree of bus $j$ in scenario $s$ at time $t$. $U^{+}$ is the upper limit of the three-phase unbalance degree. $\varepsilon$ is the risk level of exceeding the unbalance limit. $v_{j,t,s}^{avg}$ is the average three-phase voltage of bus $j$ in scenario $s$ at time $t$. Constraint (41) indicates that the risk level of exceeding the three-phase unbalance limit cannot be greater than $\varepsilon$.

*3.2.3 Model linearization*

The above model involves many nonlinear constraints, which may lead to computational difficulties. Therefore, linearization is necessary. For the circular constraints (34) and (35), the circumscribed octagon method can be employed for linearization, and the linearized formulations are given as follows:

$$\begin{cases} -S_d/\sqrt{3} \le p_{*,d,\varphi,t,s} \le S_d/\sqrt{3} \\ -S_d/\sqrt{3} \le q_{*,d,\varphi,t,s} \le S_d/\sqrt{3} \\ -S_d/\sqrt{3} \le p_{*,d,\varphi,t,s} + q_{*,d,\varphi,t,s} \le S_d/\sqrt{3} \\ -S_d/\sqrt{3} \le p_{*,d,\varphi,t,s} - q_{*,d,\varphi,t,s} \le S_d/\sqrt{3} \end{cases} \quad \forall s,t,d,\varphi,\forall * \in \{I,II\} \tag{43}$$

For the nonlinear constraints (41) and (42), the three-phase unbalance constraints are first equivalently reformulated as follows:

$$\begin{cases} v_{j,\varphi,t,s} - v_{j,t,s}^{avg} + U_{j,t,s} v_{j,t,s}^{avg} \ge 0 \\ v_{j,\varphi,t,s} - v_{j,t,s}^{avg} - U_{j,t,s} v_{j,t,s}^{avg} \le 0 \end{cases} \quad \forall s,t,\varphi,j \tag{44}$$

Then, by introducing binary variables $z_s$, the chance constraints are transformed into the following tractable formulation:

$$\begin{cases} v_{l,\varphi,t,s} - v_{l,t,s}^{avg} + U^+ v_{l,t,s}^{avg} \ge -M_1 z_s \\ v_{l,\varphi,t,s} - v_{l,t,s}^{avg} - U^+ v_{l,t,s}^{avg} \le M_1 z_s \end{cases} \quad \forall s,t,\varphi,l \tag{45}$$

$$\sum_{s\in\mathcal{S}} \pi_s z_s \le 1-\varepsilon \tag{46}$$

where $M_1$ is a sufficiently large constant to ensure that constraint holds.

*3.2.4 Compact form of the independent operation model of the DSO*

Let $\boldsymbol{y} = \left\{ S_d^{\text{SOP}}; p_{d,\varphi,t,s}; q_{d,\varphi,t,s}, \forall s,t,\varphi,d \right\}$ denote the decision variables of the DSO, then the independent operation model $\mathbf{DSO^0}$ can be formulated as follows:

$$\begin{cases} \min C^{\text{DSO}}(y) \\ \text{s.t.} \quad \boldsymbol{y} \in \mathcal{Y} = \{(29)-(33),(36)-(40),(43),(45)-(46)\} \end{cases} \tag{47}$$

## 4. DSO-FRA coordinated operation model

*4.1. Coordinated operation framework based on GNB*

In the coordinated mode, each FRA shares its respective feasible region (e.g., charging power and energy constraints for EVAs, as well as shiftable and curtailable load constraints for LAs) with the DSO. The DSO then determines the SOP deployment and operation strategies $\boldsymbol{y}$ and the FRA operation strategies $\boldsymbol{x}$ by considering both the PDN operational status and the feasible regions of FRAs, while also providing corresponding economic incentives to each FRA based on its contribution to phase balancing. To facilitate cooperation and further ensure a fair allocation of profits, this paper constructs a DSO-FRA coordination framework based on the GNB theory. Building upon the independent operation models (20) and (47), the proposed coordinated operation model is formulated as follows:

$$\max_{\boldsymbol{x},\boldsymbol{y},\boldsymbol{\kappa}} \left( \bar{C}^{\text{DSO}} - \tilde{C}^{\text{DSO}} - \sum_{u\in\mathcal{U}} \kappa_u \right)^{\alpha} \prod_{u\in\mathcal{U}} \left[ \sum_{s\in\mathcal{S}} \pi_s \left( \bar{C}_{u,s}^{\text{FRA}} - \tilde{C}_{u,s}^{\text{FRA}} \right) + \kappa_u \right]^{\beta_u} \tag{48}$$

$$\text{s.t.} \quad \bar{C}^{\text{DSO}} - \tilde{C}^{\text{DSO}} - \sum_{u\in\mathcal{U}} \kappa_u \ge 0 \tag{49}$$

$$\sum_{s\in\mathcal{S}}\pi_s\left(\bar{C}_{u,s}^{\mathrm{FRA}}-\tilde{C}_{u,s}^{\mathrm{FRA}}\right)+\kappa_u\geq 0 \quad \forall u \tag{50}$$

$$\boldsymbol{x}_{u,s}\in\mathcal{X}_{n,s}^{\mathrm{EVA}} \quad \forall s,\forall u\in\mathcal{U}\cap\mathcal{N} \tag{51}$$

$$\boldsymbol{x}_{u,s}\in\mathcal{X}_{m,s}^{\mathrm{LA}} \quad \forall s,\forall u\in\mathcal{U}\cap\mathcal{M} \tag{52}$$

$$\boldsymbol{y}\in\mathcal{Y} \tag{53}$$

where $\bar{C}^{\mathrm{DSO}}$ and $\tilde{C}^{\mathrm{DSO}}$ are the total investment and operation costs of the DSO in the independent and coordinated modes, respectively. $\bar{C}_{u,s}^{\mathrm{FRA}}$ and $\tilde{C}_{u,s}^{\mathrm{FRA}}$ are the total operation costs of FRA $u$ in scenario $s$ under the independent and coordinated modes, respectively. $\kappa_u$ is the incentive payment provided by the DSO to FRA $u$. $\alpha$ and $\beta_u$ denote the bargaining power of the DSO and FRA $u$, respectively. The objective function (48) aims to maximize the generalized Nash product of the profit increases of all participants. Constraints (49) and (50) ensure that the profit of each participant is nonnegative, which represents the individual rationality condition for participating in the coalition. Constraints (51)–(53) define the feasible regions of all participants to prevent infeasible operation.

*4.2. Bargaining power model*

In the standard Nash bargaining problem, the bargaining power parameters of all participants are set to 1, leading to identical payments to each participant. This approach fails to reflect the differences in contributions of different participants to phase balancing. To enhance the fairness of profit allocation, this paper designs the bargaining power of different participants in a differentiated manner to reflect the magnitude of their contributions to phase balancing.

First, to evaluate the flexibility contribution of different FRAs to phase balancing, this paper introduces the flexibility contribution index $\lambda$, which is defined by (54)–(56). The flexibility contribution index is calculated as a weighted sum of the cumulative energy deviation $\gamma$ and the total energy deviation $\delta$. A larger $\lambda_u$ implies that the cumulative energy profile or the total energy consumption of FRA $u$ has deviated more significantly, indicating that it has contributed greater flexibility for mitigating three-phase unbalance.

$$\lambda_u=k_1\gamma_u+k_2\delta_u \quad \forall u \tag{54}$$

$$\gamma_u=\frac{\sum_{t\in\mathcal{T}}\left|\sum_{\tau\leq t}\sum_{s\in\mathcal{S}}\pi_s\bar{x}_{u,\tau,s}-\sum_{\tau\leq t}\sum_{s\in\mathcal{S}}\pi_s\tilde{x}_{u,\tau,s}\right|\Delta t}{x_u^+T} \quad \forall u \tag{55}$$

$$\delta_u=\frac{\left|\sum_{t\in\mathcal{T}}\sum_{s\in\mathcal{S}}\pi_s\bar{x}_{u,\tau,s}-\sum_{t\in\mathcal{T}}\sum_{s\in\mathcal{S}}\pi_s\tilde{x}_{u,t,s}\right|\Delta t}{x_u^+} \quad \forall u \tag{56}$$

where $\bar{x}_{u,t,s}$ and $\tilde{x}_{u,t,s}$ are the power of FRA $u$ at time $t$ in scenario $s$ under the independent and coordinated modes, respectively.

Second, in addition to the flexibility contribution, the change in costs of participants is another factor that needs to be considered in the design of bargaining power. For instance, some FRAs may shift the majority of their loads from low-price periods to periods with slightly higher prices, thereby providing substantial flexibility with only minor additional costs. In contrast, other FRAs may need to shift a small amount of loads from low-price periods to extremely high-price periods, offering less flexibility but incurring significantly higher additional costs. For such FRAs, greater compensation is required to prevent them from abandoning the coalition. Therefore, this paper further introduces the cooperation dependence index $\sigma$, which is defined by (57)–(59). A lower cooperation dependence index indicates that a participant makes greater cost sacrifices for cooperation, implying a lower willingness to participate in the coalition.

$$\sigma^{\text{DSO}} = \frac{\bar{C}^{\text{DSO}} - \tilde{C}^{\text{DSO}}}{\bar{C}^{\text{DSO}}} \tag{57}$$

$$\sigma^{\text{FRA}} = \frac{\sum_{u\in\mathcal{U}}\sum_{s\in\mathcal{S}}\pi_s \bar{C}_{u,s}^{\text{FRA}} - \sum_{u\in\mathcal{U}}\sum_{s\in\mathcal{S}}\pi_s \tilde{C}_{u,s}^{\text{FRA}}}{\sum_{u\in\mathcal{U}}\sum_{s\in\mathcal{S}}\pi_s \bar{C}_{u,s}^{\text{FRA}}} \tag{58}$$

$$\sigma_u^{\text{FRA}} = \frac{\sum_{s\in\mathcal{S}}\pi_s \bar{C}_{u,s}^{\text{FRA}} - \sum_{s\in\mathcal{S}}\pi_s \tilde{C}_{u,s}^{\text{FRA}}}{\sum_{s\in\mathcal{S}}\pi_s \bar{C}_{u,s}^{\text{FRA}}} \quad \forall u \tag{59}$$

Finally, taking into account both the flexibility contribution index and cooperation dependence index, and adhering to the following principles: 1) the bargaining power is positively correlated with the flexibility contribution index; 2) the bargaining power is negatively correlated with the cooperation dependence index; and 3) the bargaining power parameters of all participants sum to 1 and are all positive, the final bargaining power model is formulated as follows:

$$\alpha = \frac{e^{-\sigma^{\text{DSO}}}}{e^{-\sigma^{\text{DSO}}} + e^{-\sigma^{\text{FRA}}}} \tag{60}$$

$$\beta_u = \frac{e^{-\sigma^{\text{FRA}}}}{e^{-\sigma^{\text{DSO}}} + e^{-\sigma^{\text{FRA}}}} \frac{k_3\lambda_u + k_4 e^{-\sigma_u^{\text{FRA}}}}{\sum_{u\in\mathcal{U}}\left(k_3\lambda_u + k_4 e^{-\sigma_u^{\text{FRA}}}\right)} \quad \forall u \tag{61}$$

*4.3. Complete DSO-FRA coordinated operation model*

Based on the formulations in Sections 4.1 and 4.2, the complete DSO-FRA coordinated operation model $\mathbf{MP^0}$ is formulated as follows:

$$\begin{cases} \max\limits_{x,y,\kappa} \left( \bar{C}^{\text{DSO}} - \tilde{C}^{\text{DSO}} - \sum\limits_{u\in\mathcal{U}} \kappa_u \right)^{\alpha} \prod\limits_{u\in\mathcal{U}} \left[ \sum\limits_{s\in\mathcal{S}} \pi_s \left( \bar{C}_{u,s}^{\text{FRA}} - \tilde{C}_{u,s}^{\text{FRA}} \right) + \kappa_u \right]^{\beta_u} \\ \text{s.t.} \quad (49)-(61) \end{cases} \tag{62}$$

## 5. Solution methods

*5.1. Solution method for the independent operation models*

For the independent operation models, the problem $\mathbf{DSO^0}$ is a mixed-integer linear programming (MILP) problem that can be readily solved using commercial solvers. In contrast, $\mathbf{FRA^0}$ constitutes a NEP, which is difficult to solve directly. Therefore, this paper introduces a distributed PDA to solve it.

Since the scenarios are mutually independent, a scenario-wise solution approach is adopted. First, for scenario $s$, $\mathbf{FRA^0}$ is reformulated into the following regularized form to enhance convergence:

$$\min_{\bar{x}_{u,s}} \bar{C}_{u,s}^{\text{FRA}}\left(\bar{x}_{u,s}, \bar{x}_{-u,s}\right) + \frac{\omega}{2}\left\|\bar{x}_{u,s} - \bar{x}_{u,s}^{k-1}\right\|_2^2 \quad \forall s,u \tag{63}$$

Then, at each iteration, each FRA is solved sequentially, and its operating strategy is updated according to (64) until the convergence condition is reached.

$$\bar{\boldsymbol{x}}_{u,s}^{k}=\arg\min_{\bar{\boldsymbol{x}}_{u,s}}\bar{C}_{u,s}^{\mathrm{FRA}}\left(\bar{\boldsymbol{x}}_{u,s},\bar{\boldsymbol{x}}_{-u,s}^{k-1}\right)+\frac{\omega}{2}\left\|\bar{\boldsymbol{x}}_{u,s}-\bar{\boldsymbol{x}}_{u,s}^{k-1}\right\|_2^2 \quad \forall s,u \tag{64}$$

As proven in [44], the sequence $\{\bar{\boldsymbol{x}}_{u,s}^{k}\}_{k=1}^{\infty}$ obtained from (64) is guaranteed to converge to the solution of the NEP. The computation process of the distributed PDA is summarized in Algorithm 1.

**Algorithm 1** Distributed proximal decomposition algorithm

**Initialize:** Set initial power $\bar{\boldsymbol{x}}_{u,s}^{0}$, and regularization parameter $\omega$.

- **for** $s \leftarrow 1$ **to** $S$ **do**
  - $k \leftarrow 1$
  - **for** $u \leftarrow 1$ **to** $U$ **do**
    - Solve (64)
  - **end for**
  - **while** $\left\|\bar{\boldsymbol{x}}_{s}^{k}-\bar{\boldsymbol{x}}_{s}^{k-1}\right\|_\infty \le 0.0001$ **do**
    - $k \leftarrow k+1$
    - **for** $u \leftarrow 1$ **to** $U$ **do**
      - Solve (64)
    - **end for**
  - **end while**
- **end for**
- **return** Independent operation strategy of FRAs $\bar{\boldsymbol{x}}_{s}$.

*5.2. Solution method for the coordinated operation model*

*5.2.1. Problem decomposition*

The objective of the coordinated operation model $\mathbf{MP^0}$ is formulated as a product, rendering it highly nonlinear and difficult to solve directly. Fortunately, it has been proven in [45] that this problem can be equivalently decomposed into two subproblems that can be solved sequentially, namely the social welfare maximization subproblem $\mathbf{SP^1}$ and the payment bargaining subproblem $\mathbf{SP^2}$.

Social welfare maximization subproblem $\mathbf{SP^1}$:

$$\begin{cases}\min\limits_{\boldsymbol{x},\boldsymbol{y}} C^{\mathrm{DSO}}+\sum\limits_{s\in\mathcal{S}}\sum\limits_{u\in\mathcal{U}}\pi_s C_{u,s}^{\mathrm{FRA}}\\ \text{s.t.}\quad (51)-(53)\end{cases} \tag{65}$$

Payment bargaining subproblem $\mathbf{SP^2}$:

$$\begin{cases}\max\limits_{\kappa}\left(\bar{C}^{\mathrm{DSO}}-\tilde{C}^{\mathrm{DSO}}-\sum\limits_{u\in\mathcal{U}}\kappa_u\right)^{\alpha}\prod\limits_{u\in\mathcal{U}}\left[\sum\limits_{s\in\mathcal{S}}\pi_s\left(\bar{C}_{u,s}^{\mathrm{FRA}}-\tilde{C}_{u,s}^{\mathrm{FRA}}\right)+\kappa_u\right]^{\beta_u}\\ \text{s.t.}\quad (49)-(50),(54)-(61)\end{cases} \tag{66}$$

For $\mathbf{SP^1}$, it remains a large-scale mixed-integer quadratic programming (MIQP) problem, and solving it directly poses significant computational challenges. An accelerated solution method for $\mathbf{SP^1}$ will be introduced in the next subsection. As for $\mathbf{SP^2}$, although it is still nonlinear, its closed-form analytical solution can be derived as follows:

First, by taking the logarithm of (67), it is equivalently reformulated as follows:

$$\max L = \alpha \ln\left( \bar{C}^{\text{DSO}} - \tilde{C}^{\text{DSO}} - \sum_{u\in\mathcal{U}} \kappa_u \right) + \sum_{u\in\mathcal{U}} \beta_u \ln\left( \sum_{s\in\mathcal{S}} \pi_s \left( \bar{C}^{\text{FRA}}_{u,s} - \tilde{C}^{\text{FRA}}_{u,s} \right) + \kappa_u \right) \tag{67}$$

At the optimal payments, the partial derivative of $L$ with respect to $\beta_u$ is equal to zero:

$$\frac{\partial L}{\partial \beta_u} = \frac{-\alpha}{\bar{C}^{\text{DSO}} - \tilde{C}^{\text{DSO}} - \sum_{u\in\mathcal{U}} \kappa_u} + \frac{\beta_u}{\sum_{s\in\mathcal{S}} \pi_s \left( \bar{C}^{\text{FRA}}_{u,s} - \tilde{C}^{\text{FRA}}_{u,s} \right) + \kappa_u} = 0 \tag{68}$$

Introduce auxiliary variable as follows:

$$\mu = \frac{\alpha}{\bar{C}^{\text{DSO}} - \tilde{C}^{\text{DSO}} - \sum_{u\in\mathcal{U}} \kappa_u} \tag{69}$$

Then the optimal payments can be expressed as:

$$\kappa_u = \sum_{s\in\mathcal{S}} \pi_s \left( \tilde{C}^{\text{FRA}}_{u,s} - \bar{C}^{\text{FRA}}_{u,s} \right) + \frac{\beta_u}{\mu} \tag{70}$$

In view of $\alpha + \sum_{u\in\mathcal{U}} \beta_u = 1$, the summation of (70) leads to:

$$\sum_{u\in\mathcal{U}} \kappa_u = \sum_{u\in\mathcal{U}} \sum_{s\in\mathcal{S}} \pi_s \left( \tilde{C}^{\text{FRA}}_{u,s} - \bar{C}^{\text{FRA}}_{u,s} \right) + \frac{1-\alpha}{\mu} \tag{71}$$

Finally, by combining (69), (70), (71), and eliminating $\mu$ and $\sum_{u\in\mathcal{U}} \kappa_u$, the closed-form analytical solution for the optimal payments can be expressed as follows:

$$\kappa_u = \sum_{s\in\mathcal{S}} \pi_s \left( \tilde{C}^{\text{FRA}}_{u,s} - \bar{C}^{\text{FRA}}_{u,s} \right) + \beta_u \left[ \bar{C}^{\text{DSO}} - \tilde{C}^{\text{DSO}} + \sum_{u\in\mathcal{U}} \sum_{s\in\mathcal{S}} \pi_s \left( \bar{C}^{\text{FRA}}_{u,s} - \tilde{C}^{\text{FRA}}_{u,s} \right) \right] \quad \forall u \tag{72}$$

*5.2.2. Improved bilinear Benders decomposition*

For the social welfare maximization subproblem $\mathbf{SP^1}$, it is reformulated as follows:

$$\begin{cases} \min\limits_{x,y} C^{\text{inv}} + \sum_{s\in\mathcal{S}} \pi_s \left( C^{\text{loss}}_s + \sum_{u\in\mathcal{U}} C^{\text{FRA}}_{u,s} \right) \\ \text{s.t.} \quad (51)-(53) \end{cases} \tag{73}$$

The first term of the objective function represents the investment cost of SOPs, which involves only integer variables, while the second term represents the operation cost of the PDN and FRAs, which involves only continuous variables. Therefore, following the general Benders decomposition framework [46], $\mathbf{SP^1}$ can be decomposed into an investment master problem and a series of operation subproblems, and the original problem can be solved efficiently by iteratively solving the master problem and the subproblems.

Let $\hat{S}^k_d$ and $\hat{z}^k_s$ be the fixed decision variables passed from the master problem to the subproblem $s$ at the $k$-th iteration; the operation subproblem $\mathbf{SP}^{(k)}_s$ can then be formulated as follows:

$$\begin{cases} \min f_s^k = C_s^{\text{loss}} + \sum_{u\in\mathcal{U}} C_{u,s}^{\text{FRA}} + M_2\rho & \forall s \\ \text{s.t.} \quad (29),(30)-(32),(36)-(39),(42),(51)-(52) & \\ \quad v_{l,\varphi,t,s} - v_{l,t,s}^{avg} + U^+ v_{l,t,s}^{avg} \geq -M_1 \hat{z}_s^k - \rho & \forall s,t,\varphi,l \\ \quad v_{l,\varphi,t,s} - v_{l,t,s}^{avg} - U^+ v_{l,t,s}^{avg} \leq M_1 \hat{z}_s^k + \rho & \forall s,t,\varphi,l \\ \quad S_d = \hat{S}_d^k & \forall d \\ \quad \rho \geq 0 & \end{cases} \tag{74}$$

where $M_2$ is a sufficiently large constant. It is noted that a slack variable $\rho$ is added to the three-phase unbalance constraints to ensure the feasibility of the subproblem. Then the optimality cut passed back from the subproblem to the master problem is given as:

$$\begin{cases} \xi_s \geq \left[ \hat{f}_s^k + \sum_{d\in\mathcal{D}} \lambda_d^k \left( S_d - \hat{S}_d^k \right) \right] z_s & \forall z_s = 1 \\ \xi_s \geq \left[ \hat{f}_s^k + \sum_{d\in\mathcal{D}} \lambda_d^k \left( S_d - \hat{S}_d^k \right) \right] (1 - z_s) & \forall z_s = 0 \end{cases} \tag{75}$$

where $\lambda_d^k$ is the dual multiplier of constraint $S_d = \hat{S}_d^k$, reflecting the marginal cost of SOP $d$. $\xi_s$ is the operation cost of scenario $s$. For the bilinear constraint (75), a conventional approach is to employ the McCormick linearization method; however, this inevitably introduces additional variables and constraints, thereby increasing the computational burden. Therefore, this paper develops an alternative linearization method for this constraint. An auxiliary variable $\zeta_s$ is defined as follows:

$$\zeta_s^k = \sum_{d\in D_1^k} \lambda_d^k S_d^+ \tag{76}$$

where $D_1^k = \left\{ d \in D : \lambda_d^k \geq 0 \right\}$. Thus, the linearized formulation of (75) can be expressed as follows:

$$\begin{cases} \xi_s \geq \left( \hat{f}_s^k - \sum_{d\in\mathcal{D}} \lambda_d^k \hat{S}_d^k \right) z_s + \sum_{d\in\mathcal{D}} \lambda_d^k S_d - \zeta_s^k (1 - z_s) & \forall z_s = 1 \\ \xi_s \geq \left( \hat{f}_s^k - \sum_{d\in\mathcal{D}} \lambda_d^k \hat{S}_d^k \right) (1 - z_s) + \sum_{d\in\mathcal{D}} \lambda_d^k S_d - \zeta_s^k z_s & \forall z_s = 0 \end{cases} \tag{77}$$

It can be observed that constraint (77) is activated only when certain conditions are satisfied. Taking the first constraint in (77) as an example, when $z_s = 1$, it degenerates into the standard optimality cut. When $z_s = 0$, the following relationship holds:

$$\sum_{d\in\mathcal{D}} \lambda_d^k S_d - \zeta_s^k \leq 0 \leq \xi_s \tag{78}$$

Therefore, this constraint is relaxed. Similarly, the other constraint is activated only when $z_s = 0$. Since $S_d^+$ is a known constant, this linearization method introduces neither additional constraints nor excessive relaxation.

After obtaining the optimality cut, the investment master problem $\mathbf{MP}^{(k)}$ at the $k$-th iteration is formulated as follows:

$$
\begin{cases}
\min C^{\text{inv}} + \sum_{s\in\mathcal{S}} \pi_s \xi_s \\
\text{s.t.} \quad (28),(35),(45),(76) \\
\quad \xi_s \geq \left( \hat{f}_s^l - \sum_{d\in\mathcal{D}} \lambda_d^l \hat{S}_d^l \right) z_s + \sum_{d\in\mathcal{D}} \lambda_d^l S_d - \zeta_s^l (1 - z_s) \quad \forall z_s = 1, l \leq k-1 \\
\quad \xi_s \geq \left( \hat{f}_s^l - \sum_{d\in\mathcal{D}} \lambda_d^l \hat{S}_d^l \right) (1 - z_s) + \sum_{d\in\mathcal{D}} \lambda_d^l S_d - \zeta_s^l z_s \quad \forall z_s = 0, l \leq k-1
\end{cases}
\tag{79}
$$

The computation process of the improved bilinear Benders decomposition is summarized in Algorithm 2.

**Algorithm 2** Improved bilinear Benders decomposition

```
Initialize: Set iteration counter k = 1, UB = +∞, LB = −∞.
while |UB - LB| / UB ≤ 0.0001 do
    Solve MP^(k)
    Update UB
    for s ← 1 to S do
        Solve SP_s^(k)
    end for
    Update LB
    k ← k + 1
end while
return Coordinated operation results C̃^DSO, C̃^FRA, y, x̃.
```

*5.3. Complete solution process*

Combining all the algorithms mentioned above, the solution procedure for the proposed model is illustrated in Fig. 3. First, by solving **SP**$^{1}$, the DSO's cost $\tilde{C}^{\text{DSO}}$, the FRA costs $\tilde{C}_{u,s}^{\text{FRA}}$, and the FRA operation strategies $\tilde{x}_{u,s}$ in the coordinated mode are obtained. Then, using the method described in Section 5.1, the DSO's cost $\bar{C}^{\text{DSO}}$, the FRA costs $\bar{C}_{u,s}^{\text{FRA}}$, and the FRA operation strategies $\bar{x}_{u,s}$ in the independent mode are acquired. Finally, by solving **SP**$^{2}$, the optimal payments $\kappa_u$ from the DSO to each FRA are determined.

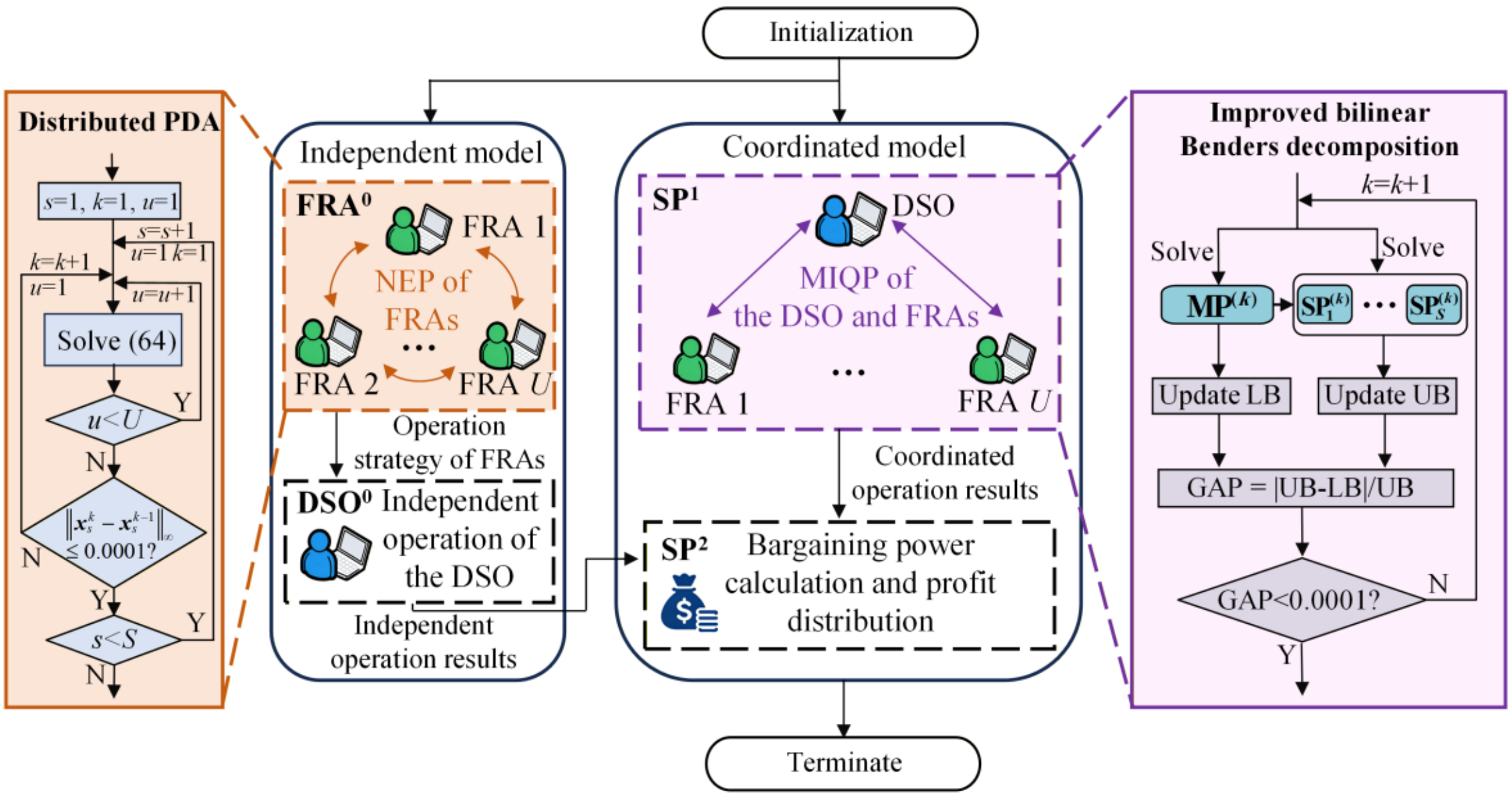

**Fig. 3.** Flow chart of the solution process.

## 6. Case studies

*6.1. Simulation parameter settings*

This section evaluates the effectiveness of the proposed method on a modified IEEE 13-bus system. The topology of the modified IEEE 13-bus system is shown in Fig. 4. The base loads and dynamic pricing coefficients are set as in [25], and all other network parameters are consistent with the standard test system [47]. The simulation period is from 13:00 to 12:00 of the next day. Each EVA manages a total of 100 EVs comprising three types, with charging powers and battery capacities of 3.8 kW/16 kWh, 3.3 kW/24 kWh, and 16.8 kW/53 kWh, respectively. The initial and expected states of charge of EVs are uniformly sampled from the ranges [0.3, 0.4] and [0.9, 1], respectively. The candidate installation locations for SOPs are buses 3–7, buses 4–8, and buses 9–10. A total of 20 typical PV scenarios are selected. The remaining parameter settings are listed in Table 1. The proposed optimization model is solved in MATLAB 2023b on a computer with a 4.8 GHz CPU and 32 GB RAM.

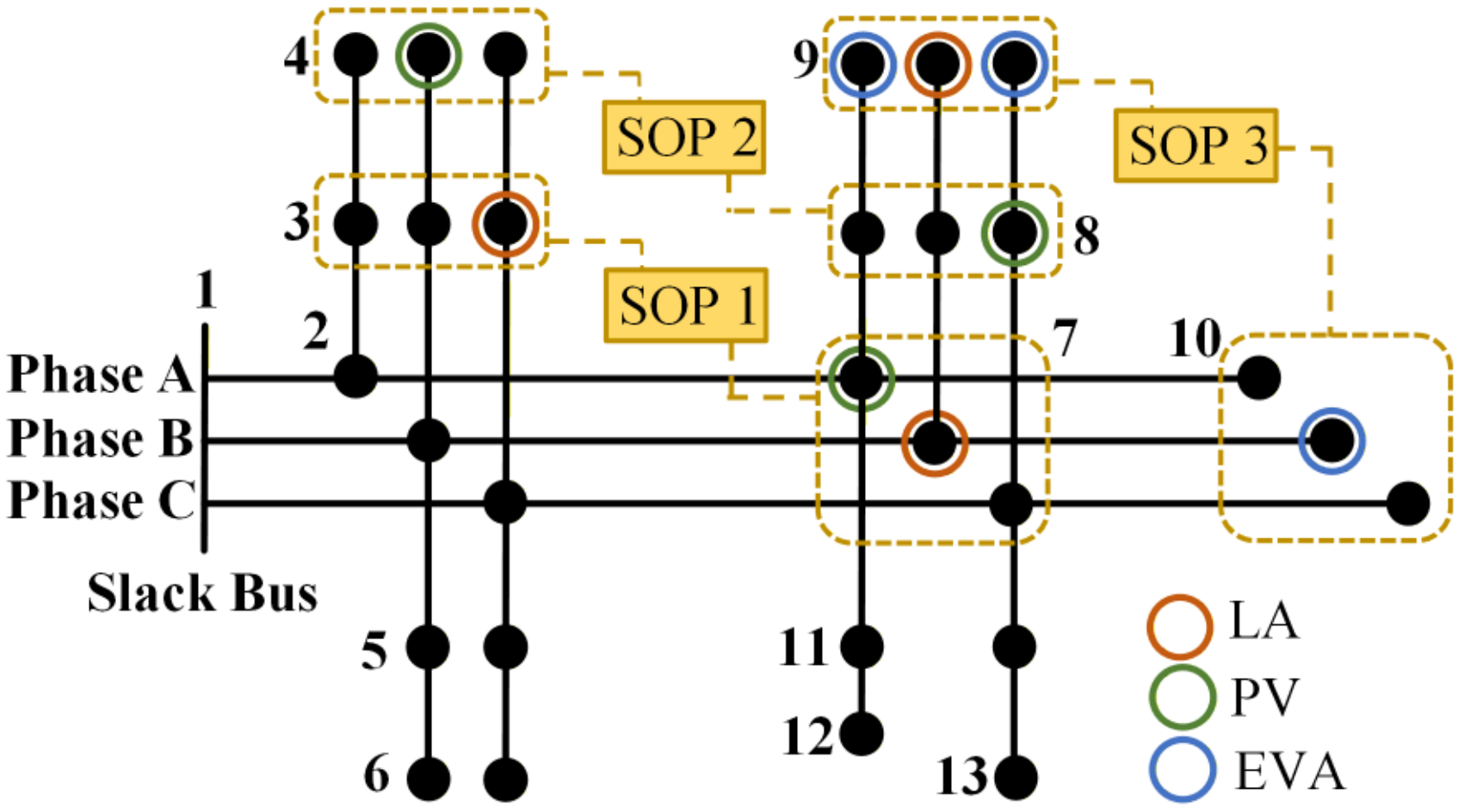

**Fig. 4.** Topology of the modified IEEE 13-bus system.

**Table 1** Simulation parameter settings

| Parameter | Value | Parameter | Value |
|---|---|---|---|
| $c_1, c_2$ | 65 \$/MWh | $r$ | 8% |
| $c_3$ | 40 \$/MWh | $y$ | 10 |
| $c_4$ | 400 \$/kVA | $S_d^+$ | 400 kVA |
| $x_n^+$ | 0.4 MW | $q^- / q^+$ | 0/400 kvar |
| $\eta$ | 0.9 | $t_i^s$ | $N$(17.6,3.42) |
| $\xi^{\text{EVA}}$, $\xi^{\text{LA}}$, $\xi^{\text{PV}}$ | tan(arccos(0.9)) | $t_i^f$ | $N$(8.0,3.22) |
| $v^- / v^+$ | 0.93/1.07 p.u. | $e_i^{\text{dev+}}$ | 10% battery capacity |
| $U^+$ | 3% | $k_1 - k_4$ | 0.5, 0.5, 0.6, 0.4 |
| $\Delta t$ | 1 h | | |

*6.2. Analysis of independent and coordinated operation modes*

In this subsection, the performance of the independent and coordinated operation modes is compared, where the risk level of unbalance violation is set to 0.

**Table 2** Cost comparison of the DSO in independent and coordinated modes

| | Independent mode (\$) | Coordinated mode (\$) | Comparison |
|---|---|---|---|
| Energy loss cost | 106.78 | 103.03 | **-3.51%** |
| Investment cost | 115.96 | 40.83 | **-64.79%** |
| Incentive payment | - | 48.22 | - |

| Total cost | 222.74 | 192.08 | **-13.76%** |
|---|---|---|---|

Table 2 compares the DSO's costs between the independent and coordinated operation modes. As shown, the energy loss cost is slightly reduced in the coordinated mode, with a decrease of 3.51%, indicating a modest improvement in network efficiency. More notably, the investment cost of SOPs drops significantly by 64.79%, demonstrating that the flexibility provided by FRAs effectively substitutes the need for extensive SOP deployment. Although the coordinated mode introduces an additional incentive payment to FRAs, the total cost of the DSO still achieves a considerable reduction of 13.76%. These results confirm that the proposed cooperative framework not only alleviates the DSO's capital expenditure on SOPs but also yields overall economic benefits, thereby validating the effectiveness of DSO-FRA coordination for unbalance mitigation.

**Table 3** Capacity comparison of SOPs in independent and coordinated modes

| | Independent mode (kVA) | Coordinated mode (kVA) | Comparison |
|---|---|---|---|
| SOP1 capacity | 40 | 30 | **-25%** |
| SOP2 capacity | 400 | 60 | **-85%** |
| SOP3 capacity | 270 | 160 | **-40.74%** |
| Total capacity | 710 | 250 | **-64.79%** |

Table 3 compares the installed capacities of SOPs between the independent and coordinated operation modes. It is evident that the total SOP capacity decreases substantially by 64.79% in the coordinated mode, indicating that the flexibility provided by FRAs effectively substitutes the need for large-scale SOP deployment. Among the individual SOPs, SOP2 experiences the most significant reduction of 85%, while SOP1 and SOP3 are also downsized by 25% and 40.74%, respectively. These results further confirm that the proposed DSO-FRA cooperative framework can significantly reduce the DSO's reliance on additional equipment investment while maintaining adequate unbalance mitigation performance.

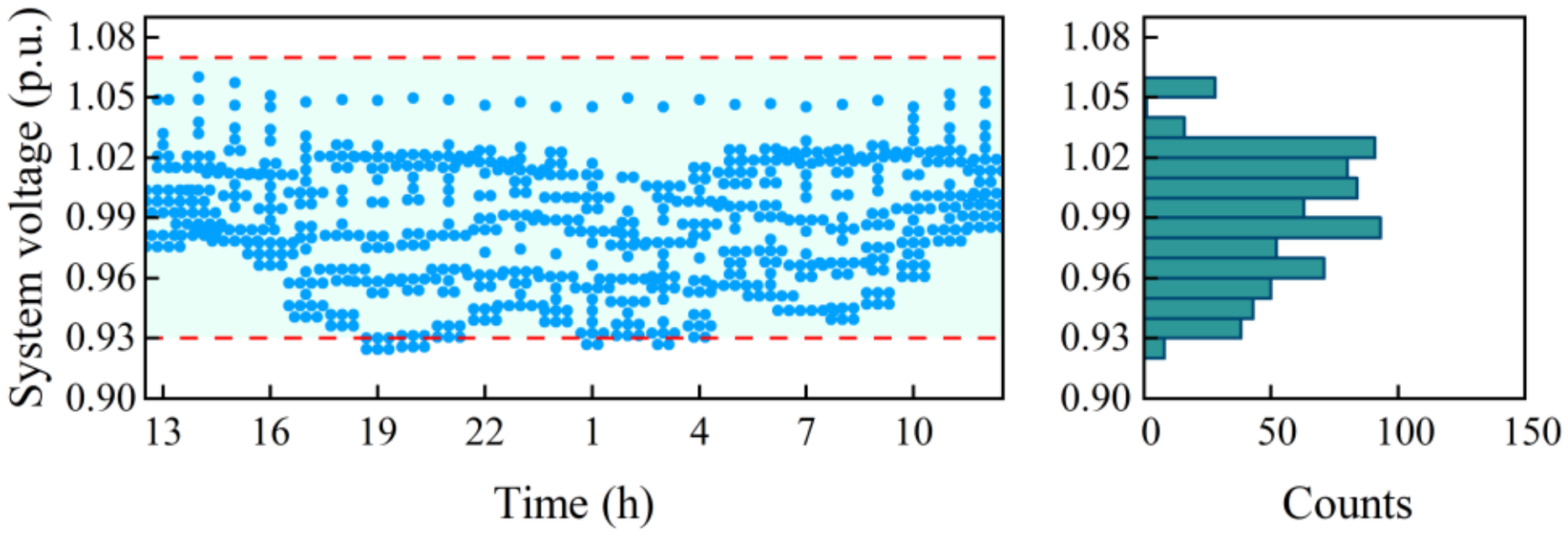

(a) Original state

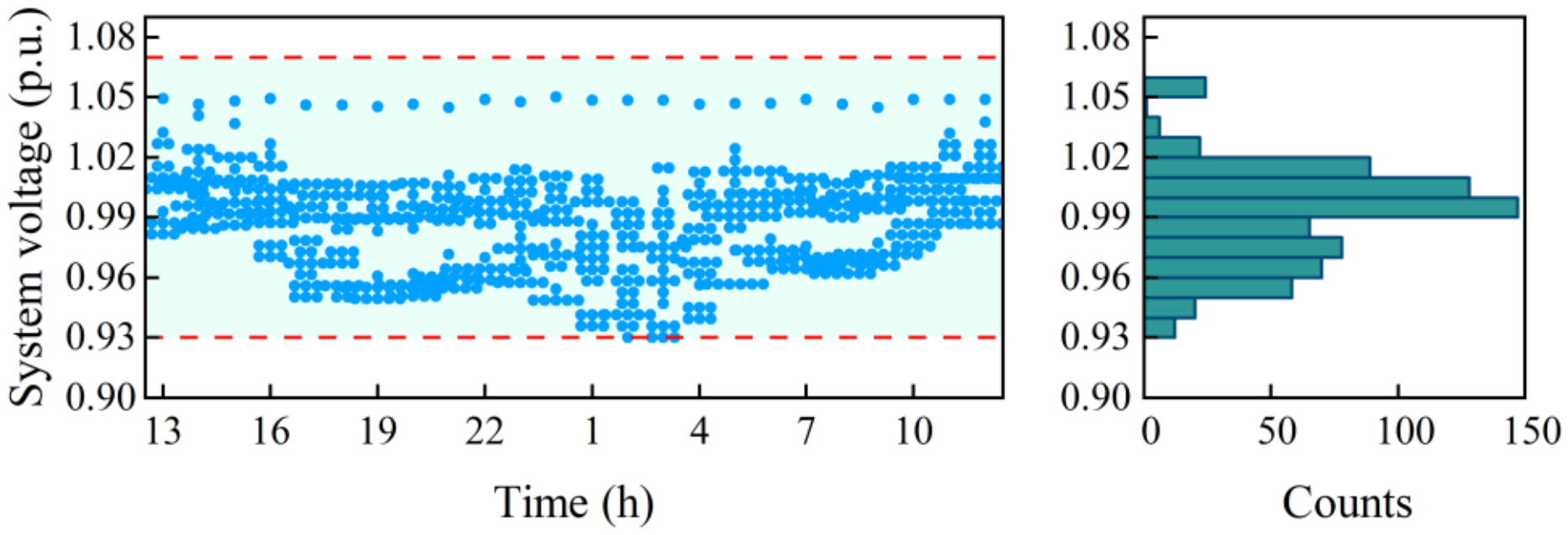

(b) Independent mode

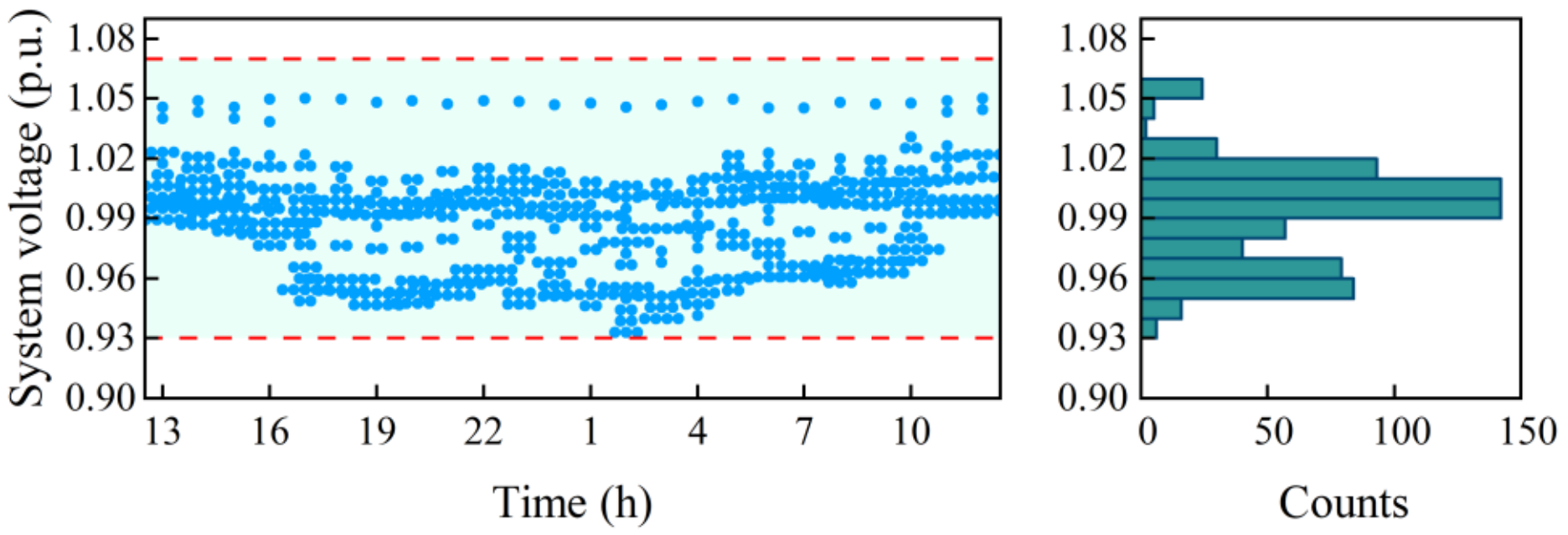

(b) Coordinated mode

**Fig. 5.** The voltage distributions of all buses in three cases.

Fig. 5 compares the voltage distributions of all buses across the original state, the independent mode, and the coordinated mode. In the original state, the voltages of certain buses violate the lower limit during certain periods. In the independent mode, through the extensive deployment of SOPs, the voltage violations are effectively eliminated. In the coordinated mode, the voltage violations are eliminated by deploying only a small number of SOPs and mobilizing the flexibility of FRAs. Notably, the voltage distribution in the coordinated mode is more concentrated around the nominal value compared to that in the independent mode, indicating that the coordinated operation of FRAs achieves a more uniform and superior voltage profile at a significantly lower cost. These results demonstrate that the proposed DSO-FRA coordination framework not only ensures voltage security but also achieves better voltage quality with reduced equipment investment.

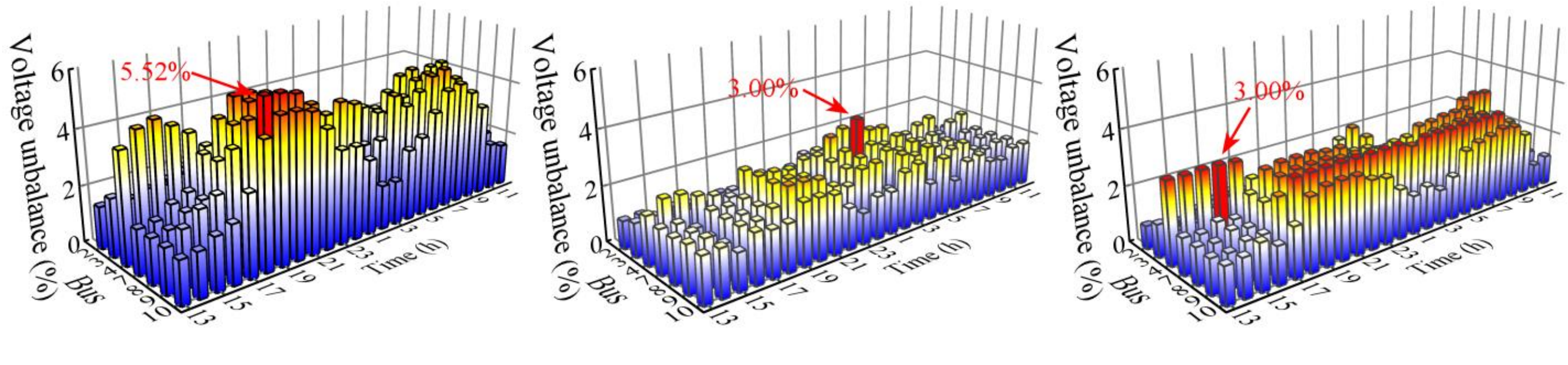

(a) Original state (b) Independent mode (c) Coordinated mode

**Fig. 6.** The voltage unbalance degrees of all three-phase buses for each case.

Fig. 6 illustrates the voltage unbalance degrees of all three-phase buses in the original state, the independent mode, and the coordinated mode. It is evident that the original state suffers from severe three-phase unbalance, with the maximum unbalance degree reaching 5.52%. In the independent mode, through the optimized deployment and operation of SOPs, the maximum unbalance degree is reduced to 3%. In the coordinated mode, by mobilizing the flexibility of FRAs, the maximum unbalance degree is also reduced to 3% with significantly lower investment costs. These results demonstrate the effectiveness of the proposed DSO-FRA coordination mechanism for three-phase unbalance mitigation.

**Table 4** Cost comparison of FRAs in independent and coordinated modes

| | FRAs | Purchase cost | Deviation cost | Incentive benefit | Total cost |
|---|---|---|---|---|---|
| Independent mode ($) | EVA 1 | 141.14 | 0 | 0 | 141.14 |
| | EVA 2 | 133.64 | 0 | 0 | 133.64 |
| | EVA 3 | 146.38 | 0 | 0 | 146.38 |
| | LA 1 | 91.81 | 0 | 0 | 91.81 |
| | LA 2 | 102.67 | 0 | 0 | 102.67 |
| | LA 3 | 103.78 | 0 | 0 | 103.75 |
| Coordinated mode ($) | EVA 1 | 121.33 | 20.93 | 9.70 | 132.56 |
| | EVA 2 | 132.45 | 1.30 | 4.67 | 129.08 |
| | EVA 3 | 146.49 | 0 | 4.50 | 141.99 |
| | LA 1 | 90.55 | 0.25 | 4.98 | 85.82 |
| | LA 2 | 104.14 | 0 | 13.74 | 90.40 |
| | LA 3 | 106.08 | 0 | 10.63 | 95.45 |

Table 4 presents a detailed cost comparison of each FRA between the independent and coordinated operation modes. It can be observed that in the independent mode, all FRAs incur only purchase costs, as they optimize their own schedules without considering system-level objectives. However, in the coordinated mode, FRAs adjust their consumption patterns to support unbalance mitigation, which leads to an increase in operational costs for most FRAs. Despite the increase in operational costs, all FRAs receive incentive payments from the DSO as compensation for their flexibility provision. As a result, the total costs of all FRAs decrease in the coordinated mode compared to the independent mode, with reductions ranging from approximately 3% to 12% across different FRAs. These results indicate that the proposed incentive mechanism effectively compensates FRAs for their flexibility contributions, ensuring that all participants benefit from the cooperation. Moreover, the different incentive payments received by FRAs reflect their differentiated flexibility contributions.

**Table 5** Comparison of bargaining power among different agents

| | Flexibility contribution index | Cooperation dependence index | Bargaining powers | Comparison |
|---|---|---|---|---|
| DSO | - | 0.3561 | 0.4105 | |
| EVA 1 | 0.8138 | -0.0079 | 0.1147 | |
| EVA 2 | 0.1242 | -0.0008 | 0.0611 | |
| EVA 3 | 0.0944 | -0.0007 | 0.0588 | |
| LA 1 | 0.3789 | 0.0110 | 0.0801 | |
| LA 2 | 1.4488 | -0.0144 | 0.1640 | |
| LA 3 | 0.7549 | -0.0225 | 0.1109 | |

Table 5 presents the flexibility contribution index, cooperation dependence index, and the resulting bargaining power for each participant. Among the FRAs, EVA 1, LA 2, and LA 3 exhibit relatively higher bargaining power, while EVA 2, EVA 3, and LA 1 show comparatively lower values. This variation primarily stems from the differentiated flexibility contributions of each FRA: those providing greater flexibility adjustments are assigned higher bargaining power. Additionally, the cooperation dependence index values further influence the bargaining power, where a positive cooperation dependence index indicates that a participant benefits from cooperation, thereby weakening its bargaining position, whereas a negative cooperation dependence index suggests that the participant sacrifices individual benefits for cooperation, thus strengthening its bargaining power. Overall, the bargaining power distribution reflects a reasonable and fair allocation principle, ensuring that participants with greater flexibility contributions and greater sacrifices made for cooperation receive correspondingly higher negotiation positions, which is essential for incentivizing active participation and sustaining the cooperative framework.

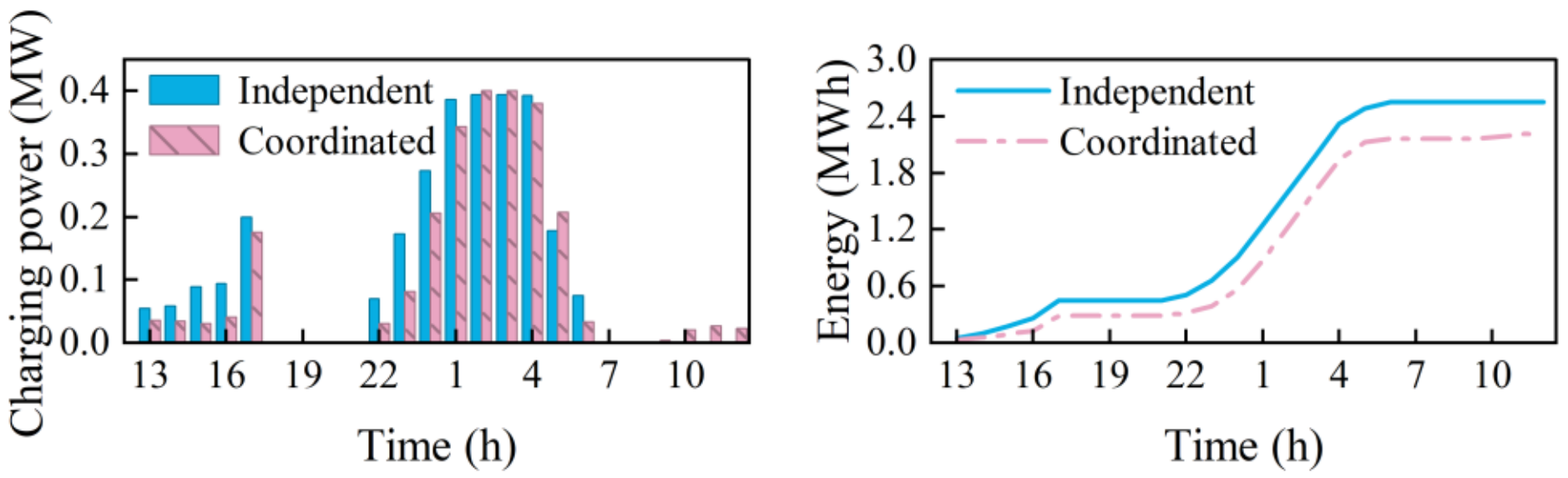

(a) EVA 1

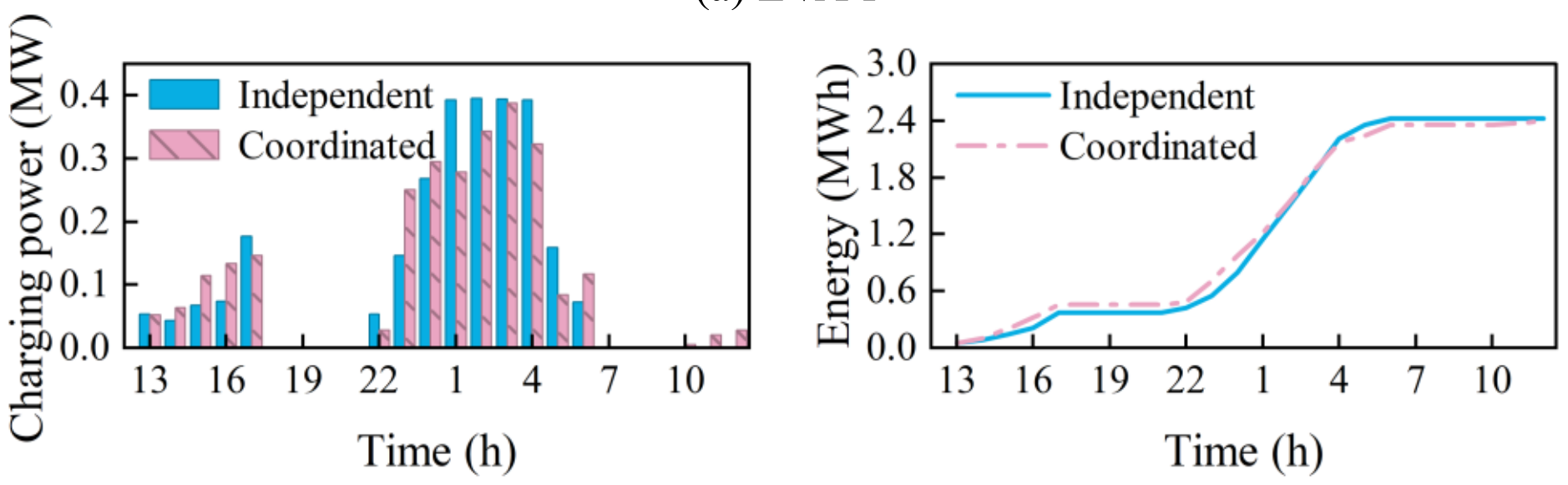

(b) EVA 2

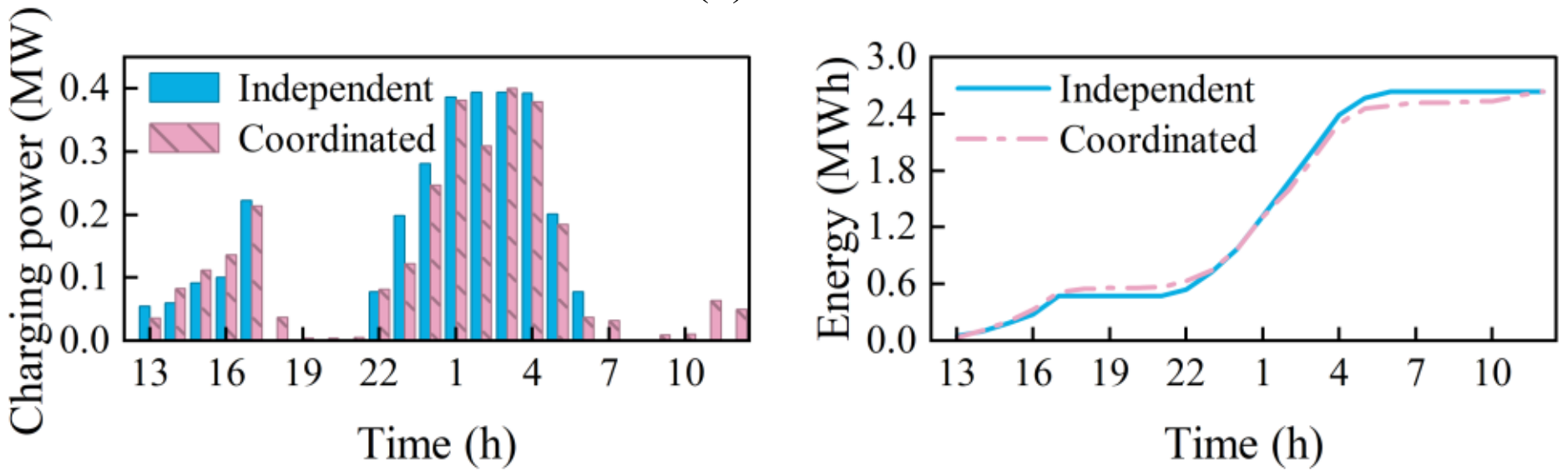

(c) EVA 3

**Fig. 7.** Charging power and energy profiles of EVAs in independent and coordinated modes.

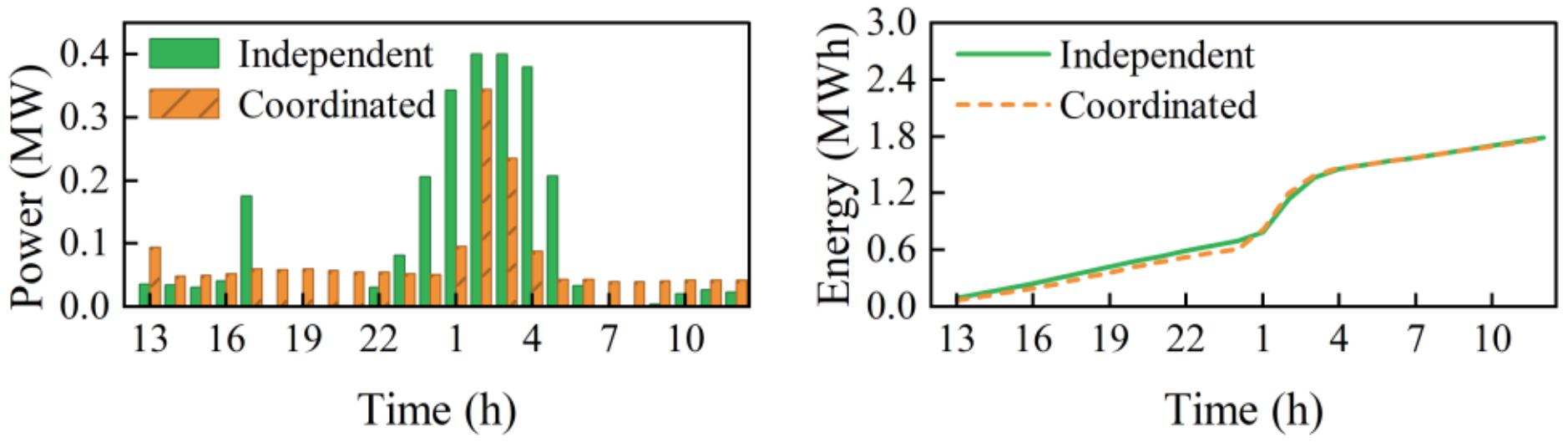

(a) LA 1

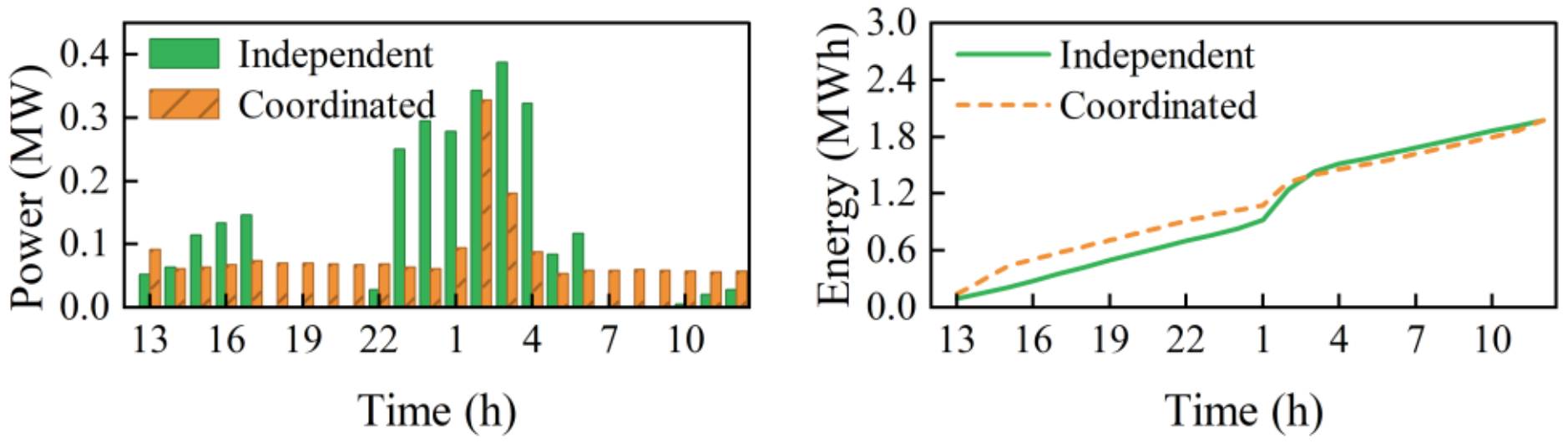

(b) LA 2

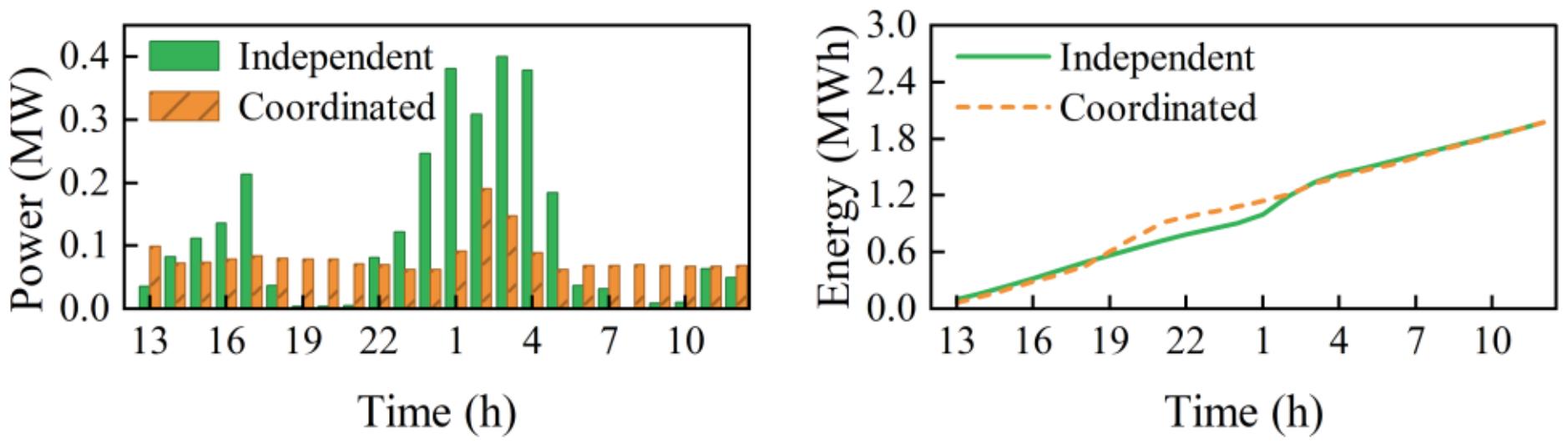

(c) LA 3

**Fig. 8.** Charging power and energy profiles of EVAs in independent and coordinated modes.

Fig. 7 and Fig. 8 present the operational profiles of FRAs in the independent and coordinated modes. As shown in Fig. 7, in the independent mode, the EVAs tend to concentrate their charging activities during low-price periods to minimize charging costs, resulting in pronounced peak charging powers. In contrast, in the coordinated mode, the charging powers of EVAs are redistributed across time slots, with EVAs shifting their charging demand away from the original peak periods to support three-phase unbalance mitigation. In addition, the accumulated energy profiles of EVAs deviate from those in the independent mode, indicating that EVAs actively provide temporal and energy flexibility to the system. Similarly, Fig. 8 illustrates that in the independent mode, the power consumption of LAs is relatively concentrated, whereas in the coordinated mode, the power profiles are dispersed across different time periods. Meanwhile, the accumulated energy trajectories of LAs also exhibit notable deviations from the independent mode, further demonstrating that LAs contribute to flexibility provision. Collectively, these results demonstrate that the proposed coordination framework successfully mobilizes the flexibility potential of both EVAs and LAs.

*6.3. Analysis of the impact of chance constraints*

In this subsection, the impact of chance constraints on the independent and coordinated operation modes is analyzed by varying the risk level of unbalance violation.

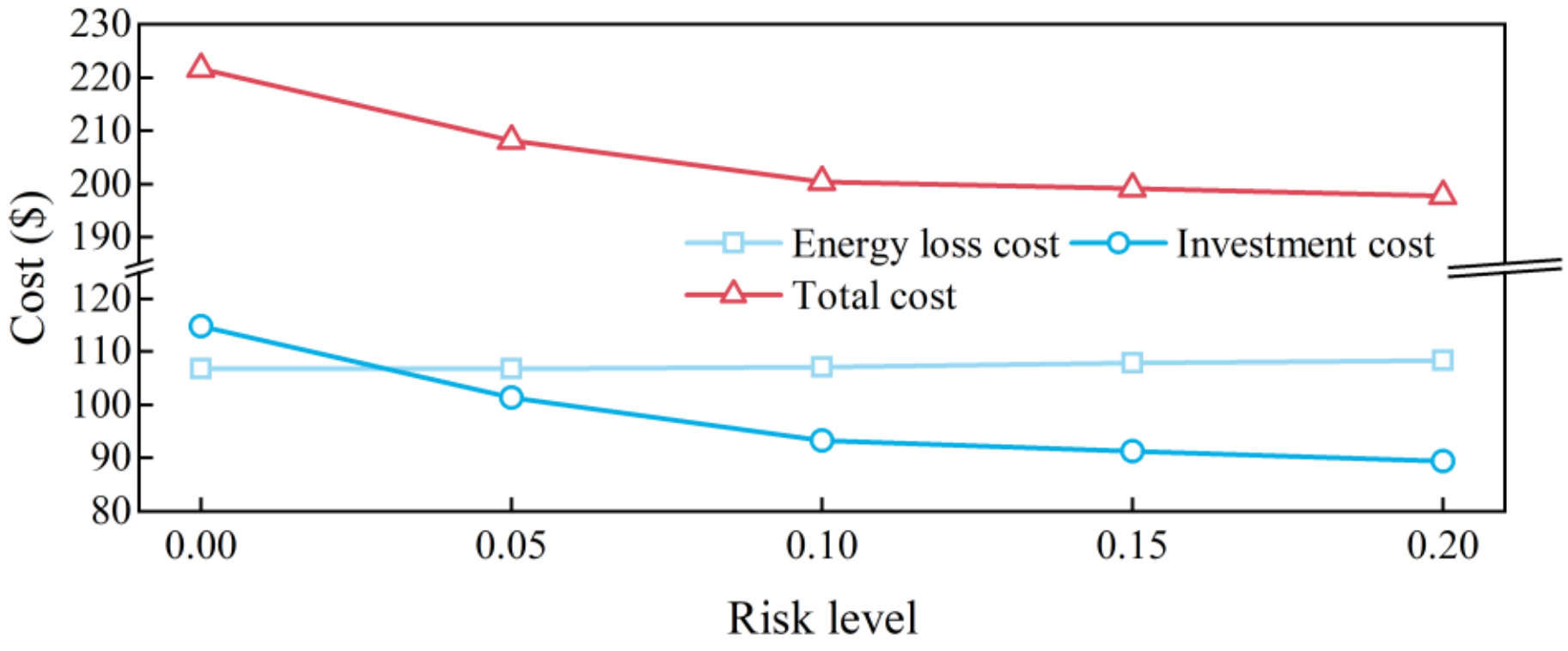

(b) Coordinated operation mode.

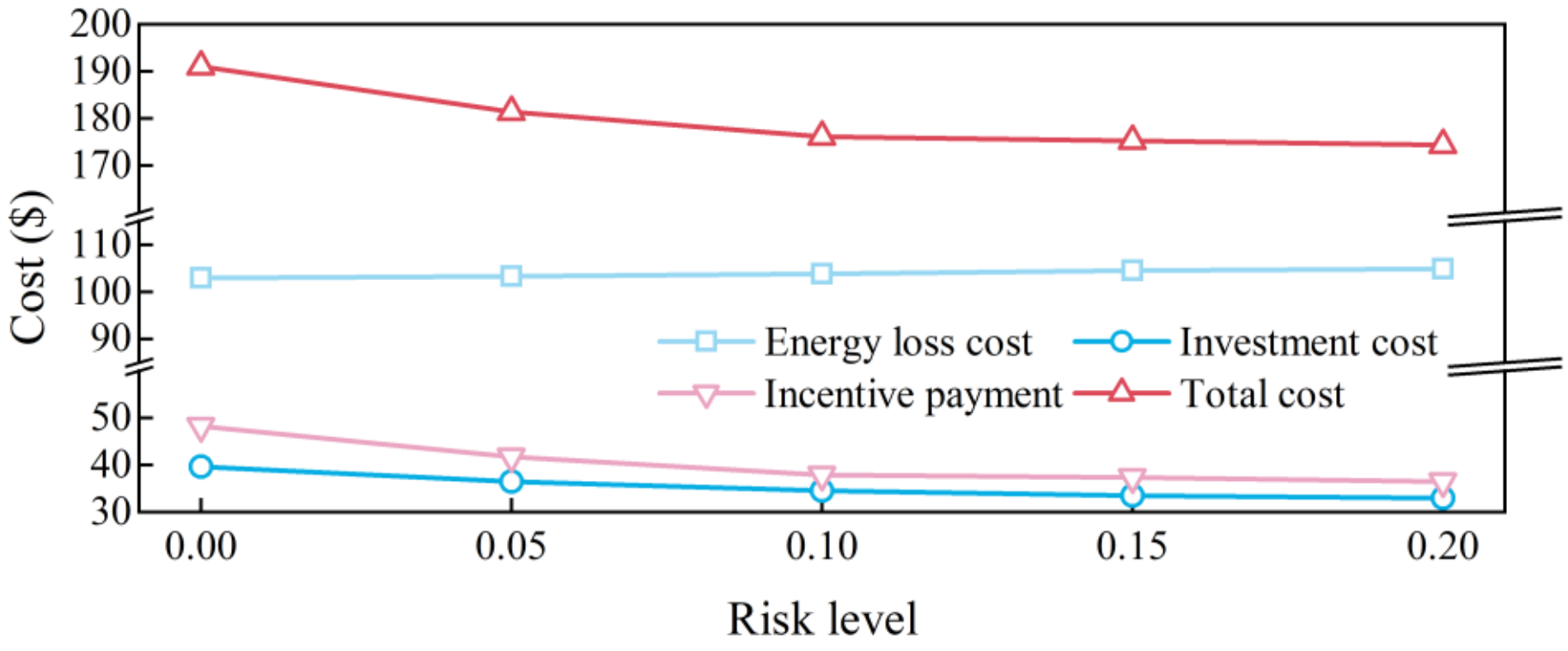

(a) Independent operation mode.

**Fig. 9.** The cost of the DSO with different risk levels in independent and coordinated modes.

Fig. 9 presents the impact of varying allowable risk levels of unbalance violation on the DSO's costs in both independent and coordinated operation modes. It can be observed that as the risk level increases, the energy loss cost exhibits a slight upward trend in both modes, indicating that relaxing the unbalance constraints tends to degrade network efficiency. In contrast, the investment cost of SOPs decreases significantly as the risk level increases in both modes, suggesting that allowing a higher probability of unbalance violation effectively reduces the DSO's need for additional equipment investment. In the coordinated mode, the incentive payments to FRAs also show a decreasing trend, implying that less stringent unbalance requirements reduce the flexibility that FRAs need to provide. Consequently, the total cost of the DSO declines steadily with increasing risk levels in both modes. Overall, the results demonstrate that the chance-constrained formulation provides a flexible trade-off between operational security and economic efficiency, enabling the DSO to achieve cost savings by accepting a certain level of unbalance violation risk, while the coordinated mode consistently outperforms the independent mode across all risk levels in terms of total cost reduction.

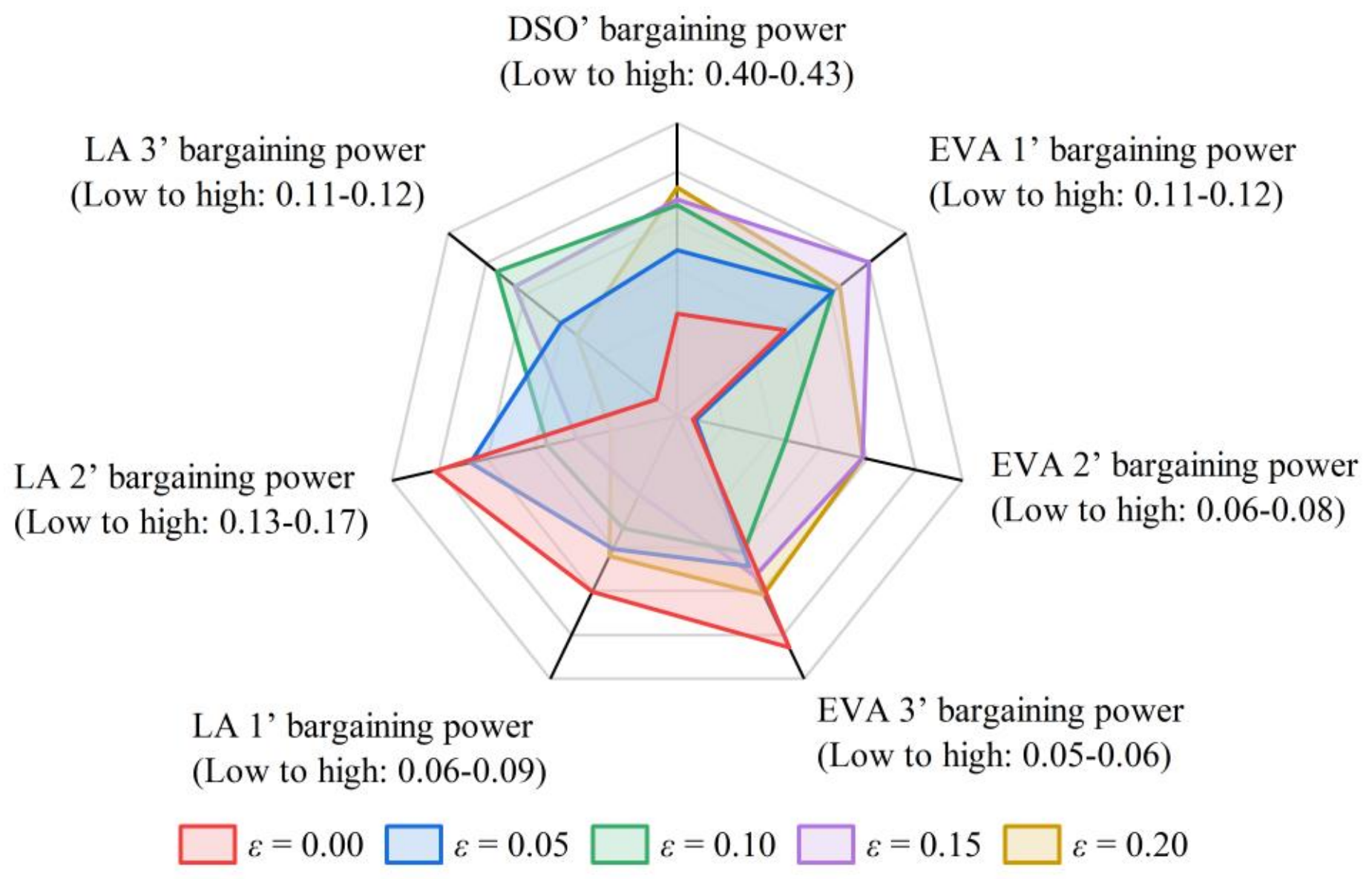

**Fig. 10.** Bargaining powers of FRAs with different risk levels.

Fig. 10 illustrates the evolution of bargaining power for each participant under different allowable risk levels of unbalance violation. It can be observed that as the risk level increases, the bargaining power of the DSO exhibits a consistent upward trend, indicating that the DSO gains a stronger negotiation position when unbalance constraints are relaxed. This is primarily because a higher tolerance for unbalance violation reduces the DSO's reliance on FRAs for flexibility provision, thereby weakening the FRAs' leverage in the bargaining process. Consequently, the flexibility contributions and the associated bargaining power of FRAs generally exhibit a decreasing trend as the risk level increases. Overall, the results suggest that the stringency of unbalance constraints not only affects the economic costs but also influences the relative bargaining positions of the participants, which is a critical factor for designing fair and sustainable incentive mechanisms under uncertain operating conditions.

### *6.4. Analysis of computational performance*

In this subsection, the computational performance of the proposed improved bilinear Benders decomposition algorithm is analyzed. First, the social welfare maximization subproblem $\mathbf{SP}^1$ is solved using both the commercial solver GUROBI and the proposed improved bilinear Benders decomposition algorithm, respectively, to verify the computational efficiency of the proposed method. During the solving process, an optimality gap of 0.01% is set for both methods. Table 6 presents the solution times of the two methods under different risk levels.

**Table 6** Comparison of different solution methods

| Risk level | Solution time (s) | | Comparison |
|---|---|---|---|
| | GUROBI | Bilinear Benders | |
| 0.00 | 73.51 | 72.06 | **-1.97%** |
| 0.05 | 287.55 | 178.29 | **-38.00%** |
| 0.10 | 464.19 | 128.24 | **-72.37%** |
| 0.15 | 731.28 | 122.06 | **-83.31%** |
| 0.20 | 930.19 | 96.49 | **-89.63%** |

It can be observed that for different risk levels, the solution time of the proposed improved bilinear Benders decomposition algorithm is consistently lower than that of GUROBI, demonstrating the superiority of the proposed method in terms of computational efficiency. Moreover, as the risk level increases, the solution time of GUROBI increases significantly, while that of the proposed method remains relatively stable. Fig. 11 illustrates the iterative process of the upper and lower bounds of the proposed method under different cases. It can be seen that the proposed method converges rapidly, with most cases reaching the convergence condition within 20 iterations, which is sufficient for practical computational requirements.

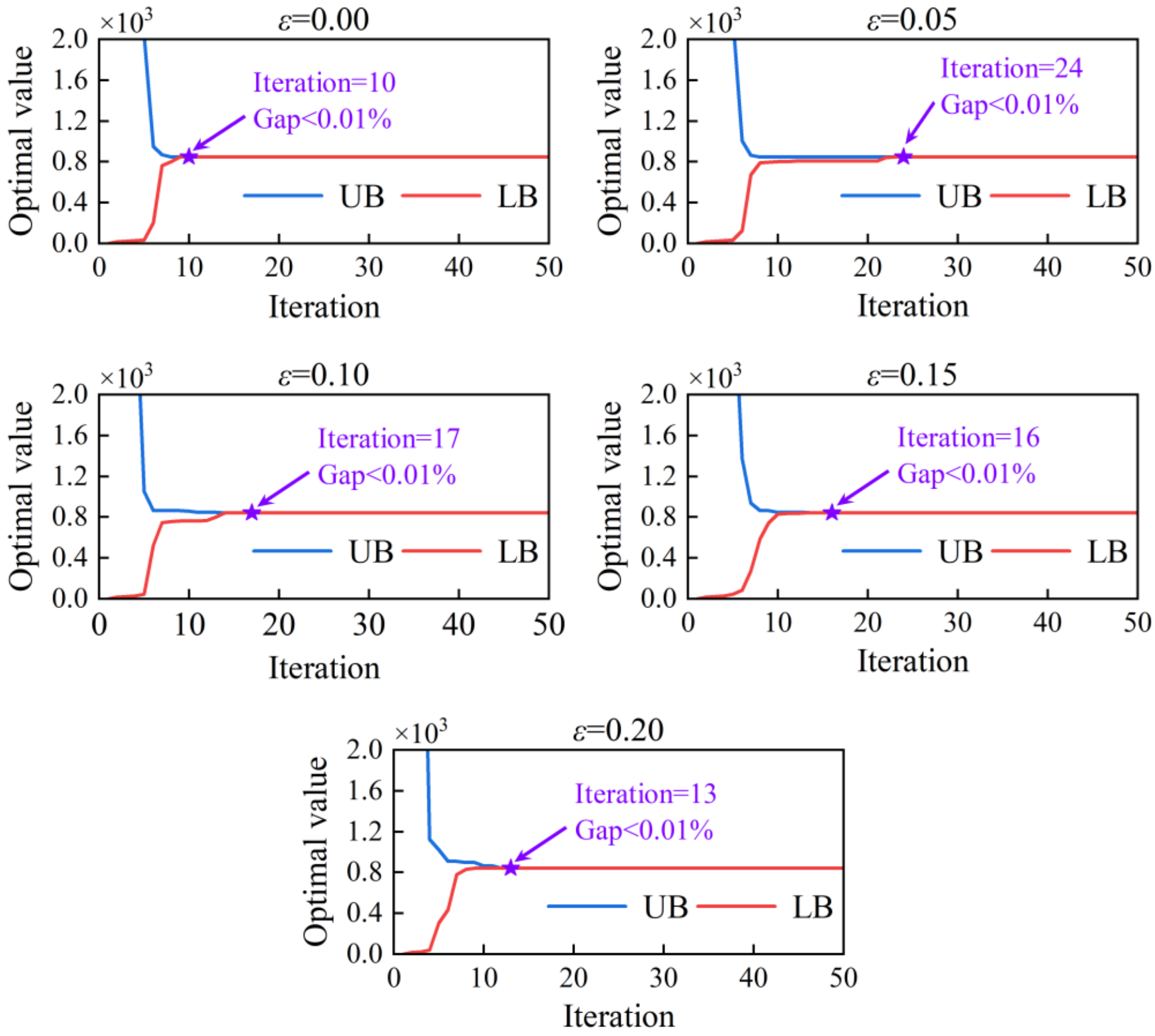

**Fig. 11.** Convergence curves of the bilinear Benders algorithm for each case.

To further validate the scalability of the proposed solution method, **SP¹** is re-solved on the basis of an increased number of scenarios, and the results are presented in Table 7. It can be observed that the average solution time per iteration exhibits an approximately linear relationship with the number of scenarios, while the number of iterations remains relatively stable. Consequently, the total computation time shows a linear increasing trend as the number of scenarios grows, indicating that the proposed algorithm possesses favorable convergence and scalability.

**Table 7** Comparison of computational performance under different numbers of scenarios

| Risk level | Number of scenarios | | | | | | | | |
|---|---|---|---|---|---|---|---|---|---|
| | 20 | | | 30 | | | 40 | | |
| | Total time (s) | Iter. | Time per Iter. (s) | Total time (s) | Iter. | Time per Iter. (s) | Total time (s) | Iter. | Time per Iter. (s) |
| 0.00 | 72.06 | 10 | 7.21 | 137.35 | 12 | 11.45 | 175.94 | 12 | 14.66 |
| 0.05 | 178.29 | 24 | 7.43 | 306.22 | 27 | 11.34 | 398.83 | 26 | 15.34 |
| 0.10 | 128.24 | 17 | 7.54 | 181.87 | 16 | 11.37 | 273.31 | 18 | 15.18 |
| 0.15 | 122.06 | 16 | 7.63 | 192.52 | 17 | 11.32 | 237.95 | 15 | 15.86 |
| 0.20 | 96.49 | 13 | 7.42 | 182.25 | 16 | 11.39 | 251.36 | 15 | 16.76 |

## 7. Conclusion

This paper proposes a GNB-based chance-constrained coordinated operation model for the DSO and FRAs to mitigate three-phase unbalance in PDNs. Numerical results on the modified IEEE 13-bus system demonstrate that the proposed method effectively mitigates three-phase unbalance, reducing the maximum unbalance degree from 5.52% to 3%, while significantly reducing the DSO's total cost by approximately 13.76% and the SOP investment cost by 64.79%. The operational costs of all FRAs are also reduced, with total cost decreases ranging from 3% to 12% across different FRAs. The cooperative framework also ensures that the surplus gains of each participant are determined by their individual contributions to unbalance mitigation, thereby confirming the fairness of the proposed incentive mechanism. Moreover, by relaxing the three-phase unbalance constraints through the chance-constrained formulation, the DSO can effectively reduce its unbalance mitigation costs and simultaneously improve its bargaining position in the negotiation with FRAs. Furthermore, the proposed bilinear Benders decomposition algorithm reduces the computation time by over 75%, significantly improving the solution efficiency compared to direct solving with commercial solvers.

Future work could extend the proposed method by considering additional uncertainty sources, such as load variability and EV charging behavior to further improve the practical applicability of the proposed approach.

## Acknowledgements

This work was financially supported by the National Natural Science Foundation of China under Grant 52607198.

## References

[1] Li T, Zhang M, Hu Z, Wang X, Zhou Y, Yan M. Distributed cooperative scheduling for distribution network and smart charging hubs driven by unbalanced distribution locational marginal price. Appl Energy 2025;401:126668.

[2] Kong W, Ma K, Fang L, Wei R, Li F. Cost-benefit analysis of phase balancing solution for data-scarce LV networks by cluster-wise Gaussian process regression. IEEE Trans Power Syst 2020;35:3170-3180.

[3] Ma K, Fang L, Kong W. Review of distribution network phase unbalance: Scale, causes, consequences, solutions, and future research directions. CSEE J Power Energy Syst 2020;6:479-488.

[4] Li P, Ji H, Wang C, Zhao J, Song G, Ding F, et al. Optimal operation of soft open points in active distribution networks under three-phase unbalanced conditions. IEEE Trans Smart Grid 2019;10:380-391.

[5] De Abreu JPG, Emanuel AE. Induction motor thermal aging caused by voltage distortion and imbalance: loss of useful life and its estimated cost. IEEE Trans Ind Appl 2002;38:12-20.

[6] Ma K, Li R, Li F. Quantification of additional asset reinforcement cost from 3-phase imbalance. IEEE Trans Power Syst 2016;31:2885-2891.

[7] Chang Y, Zhao M, Kocar I. The impact of DFIG control schemes on the negative-sequence based differential protection. Electr Power Syst Res 2022;211:108564.

[8] Woll RF. Effect of unbalanced voltage on the operation of polyphase induction motors. IEEE Trans Ind Appl 1975;IA-11:38-42.

[9] Liu S, Jin R, Qiu H, Cui X, Lin Z, Lian Z, et al. Practical method for mitigating three-phase unbalance based on data-driven user phase identification. IEEE Trans Power Syst 2020;35:1653-1656.
[10] Tian S, Jia Q, Cui Y, Xue S, Yu H, Liu W. Multi-objective collaborative optimization of VDAPFs and SVGs allocation considering MFGCIs contribution for voltage partitioning mitigation in distribution networks. Electr Power Syst Res 2022;207:107830.
[11] Gupta RK, Molzahn DK. Improving fairness in photovoltaic curtailment via feedback-driven daily topology reconfiguration in power distribution networks. Appl Energy 2025;400:126543.
[12] Guo J, Zhou Q, Huang J, Gu C. A chance-constrained planning method of phase switch devices for three-phase unbalance mitigation in low-voltage distribution networks with PV uncertainty. Electr Power Syst Res 2026;259:113221.
[13] Liu B, Meng K, Dong ZY, Wong PKC, Ting T. Unbalance mitigation via phase-switching device and static var compensator in low-voltage distribution network. IEEE Trans Power Syst 2020;35:4856-4869.
[14] Shi H, Xiong H, Gan W, Lin Y, Guo C. Fully distributed planning method for coordinated distribution and urban transportation networks considering three-phase unbalance mitigation. Appl Energy 2025;377:124449.
[15] Huang J, Wang Y, Gu C, Cheng L, Fan J, Wang X. DLMP-based congestion management model for power distribution network considering network loss and EV charging demand uncertainty. IEEE Trans Smart Grid 2026;17:1414-1429.
[16] Lepolesa LJ, Adetunji KE, Ouahada K, Liu Z, Cheng L. Dynamic electric vehicle charging pricing for load balancing in power distribution networks based on collaborative DDPG agents. IEEE Trans Smart Grid 2026;17:3342-3353.
[17] Estebsari A, Mazzarino PR, Bottaccioli L, Patti E. IoT-enabled real-time management of smart grids with demand response aggregators. IEEE Trans Ind Appl 2022;58:102-112.
[18] Silveira Junior JR, Conrado BRP, Alonso AMD, Brandao DI. Interoperability of single-controllable clusters: aggregate response of low-voltage microgrids. Appl Energy 2023;340:121042.
[19] Sumon Rashid THM, Hossain MA, Kabir H, Refat KA, Islam MN, Khan R, et al. A data-driven framework for seasonal PV forecasting and coordinated EV charging to improve hosting capacity in unbalanced distribution grids. Energy Convers Manag X 2026;30:101708.
[20] Saleh SV, Latify MA. Coordinated unbalance compensation and harmonic mitigation in the secondary distribution network through EVs participation. IEEE Trans Smart Grid 2024;15:4450-4462.
[21] Zhou S, Han Y, Zalhaf AS, Lehtonen M, Darwish MMF, Mahmoud K. A novel stochastic multistage dispatching model of hybrid battery-electric vehicle-supercapacitor storage system to minimize three-phase unbalance. Energy 2024;296:131174.
[22] Zhang M, Gong J, Wang X, Wu Q, Zhou B, Li F, et al. Coordinated scheduling of DER aggregators driven by unbalanced distribution LMP. IEEE Trans Smart Grid 2026;17:936-951.
[23] Tiwari A, Jha BK, Pindoriya NM. Incentive-based demand response program with phase unbalance mitigation: A bilevel approach. Sustain Energy Grids Netw 2025;42:101671.
[24] Chen S, Guo Z, Yang Z, Xu Y, Cheng RS. A game theoretic approach to phase balancing by plug-in electric vehicles in the smart grid. IEEE Trans Power Syst 2020;35:2232-2244.
[25] Huang J, Wang X, Wang Y, Shao C, Chen G, Wang P. A game-theoretic approach for electric vehicle aggregators participating in phase balancing considering network topology. IEEE Trans Smart Grid 2024;15:743-756.
[26] Chen Y, Park B, Kou X, Hu M, Dong J, Li F, et al. A comparison study on trading behavior and profit distribution in local energy transaction games. Appl Energy 2020;280:115941.
[27] Alizadeh A, Esfahani M, Dinar F, Kamwa I, Moeini A, Mohseni-Bonab SM, et al. A cooperative transactive multi-carrier energy control mechanism with P2P energy + reserve trading using Nash bargaining game theory under renewables uncertainty. Appl Energy 2024;353:122162.
[28] Meng Y, Ma G, Ye Y, Yao Y, Li W, Li T. Design of P2P trading mechanism for multi-energy prosumers based on generalized Nash bargaining in GCT-CET market. Appl Energy 2024;371:123640.
[29] Zhang Q, Mu Y, Jia H, Jiang Q, Yu X, et al. Multistage optimized planning for the alliance expansion of multiple island microgrids with shared hybrid energy storage systems: A Nash bargaining-based game approach. Appl Energy 2026;406:127280.
[30] Xi Y, Zuo L, Chen M, Zhang J, Cheng L, Lestas I. Energy-carbon pricing-guided collaborative optimization for local integrated energy communities with multi-vector electrification consumers. Appl Energy 2025;392:125946.
[31] Cui S, Xu S, Fang J, Ai X, Wen J. A novel stable grand coalition for transactive multi-energy management in an integrated energy system. Appl Energy 2025;394:126155.
[32] Dar MR, Ganguly S. Voltage regulation and loss minimization of active distribution networks with uncertainties using chance-constrained model predictive control. IEEE Trans Power Syst 2025;40:2737-2749.
[33] Li X, Zhang L, Wang R, Sun B, Xie W. Two-stage robust optimization model for capacity configuration of biogas-solar-wind integrated energy system. IEEE Trans Ind Appl 2023;59:662-675.
[34] Chen B, Liu T, Liu X, He C, Nan L, Wu L, et al. A Wasserstein distance-based distributionally robust chance-

constrained clustered generation expansion planning considering flexible resource investments. IEEE Trans Power Syst 2023;38:5635-5647.
[35] Yang Z, Hu J, Ai X, Wu J, Yang G. Transactive energy supported economic operation for multi-energy complementary microgrids. IEEE Trans Smart Grid 2021;12:4-17.
[36] Xu J, Wang J, Liao S, Sun Y, Ke D, Li X, et al. Stochastic multi-objective optimization of photovoltaics integrated three-phase distribution network based on dynamic scenarios. Appl Energy 2018;231:985-996.
[37] Silva JAA, López JC, Arias NB, Rider MJ, Da Silva LCP. An optimal stochastic energy management system for resilient microgrids. Appl Energy 2021;300:117435.
[38] Song S, Xiong H, Lin Y, Huang M, Wei Z, Fang Z. Robust three-phase state estimation for PV-integrated unbalanced distribution systems. Appl Energy 2022;322:119427.
[39] Liu F, Zeng J, Liu J, Liu Z, Zhang Q, Lu Y, et al. Distributionally robust optimization of power grid along mountain railway considering three-phase unbalance of traction load. Int J Electr Power Energy Syst 2026;178:111969.
[40] Gu C, Ruan J, Yang X, Huang J, Qiu Y, Wang J, et al. Advancing climate-adaptive operation and planning for renewable-rich energy systems. Renew Sustain Energy Rev 2026;228:117336.
[41] Gu C, Liu Y, Wang J, Li Q, Wu L. Carbon-oriented planning of distributed generation and energy storage assets in power distribution network with hydrogen-based microgrids. IEEE Trans Sustain Energy 2023;14:790-802.
[42] Huang J, Shao C, Gu C, Zhang S, Zhou Q, Leng M, et al. Facilitating unbalance mitigation in three-phase distribution network: a DSO-EVA coordination model using GNB-based cooperative framework. IEEE Trans Smart Grid 2025;16:5246-5261.
[43] Bernstein A, Wang C, Dall'Anese E, Le Boudec J-Y, Zhao C. Load flow in multiphase distribution networks: Existence, uniqueness, non-singularity and linear models. IEEE Trans Power Syst 2018;33:5832-5843.
[44] Scutari G, Palomar DP, Facchinei F, Pang J-S. Monotone games for cognitive radio systems. In: Johansson R, Rantzer A, eds. Distributed Decision Making and Control. London: Springer; 2012:83-112.
[45] Wang X, Zhao H, Lu H, Zhang Y, Wang Y, Wang J. Decentralized coordinated operation model of VPP and P2H systems based on stochastic-bargaining game considering multiple uncertainties and carbon cost. Appl Energy 2022;312:118750.
[46] Zhang C, Du M, Hu Z, Jiang X, Liu Y, Cheng L. Benders decomposition for the charging station location-routing problem under stochastic energy consumption: a two-stage model with a space-electricity network perspective. Appl Energy 2026;416:127926.
[47] IEEE PES test feeders. https://cmte.ieee.org/pes-testfeeders/resources/; 2026 [accessed 10 September 2026].